\documentclass[a4paper,10pt]{article}
\usepackage[english]{babel}
\usepackage{color}
\usepackage{amsmath}
\usepackage{amssymb}
\usepackage{amsfonts}
\usepackage{graphicx}
\usepackage{subfigure}
\usepackage{caption}
\usepackage[a4paper=true,pagebackref=true]{hyperref}

\hypersetup{
    colorlinks = true,
    linkcolor = red,
    anchorcolor = red,
    citecolor = blue,
    filecolor = red,
}

\definecolor{vertclair}{rgb}{0.7,0.9,0.8}
\definecolor{vertclair}{rgb}{0.7,0.9,0.8}
\definecolor{vertfonce}{rgb}{0.3,0.8,0.4}
\definecolor{vertfonce2}{rgb}{0.4,0.7,0.3}
\definecolor{verttresfonce}{rgb}{0.3,0.4,0.3}

\definecolor{bleuclair}{rgb}{0.945,0.941,0.937}
\definecolor{bleufonce}{rgb}{0.9,0.9,0.9}

\begin{document}
\newtheorem{Remark}{Remark}
\newtheorem{Lemma}{Lemma}
\newtheorem{Theorem}{Theorem}
\newtheorem{Proposition}{Proposition}
%
%\maketitle
\numberwithin{equation}{section}
\numberwithin{Lemma}{section}
\numberwithin{Theorem}{section}
\numberwithin{Proposition}{section}
\numberwithin{Remark}{section}
%\numberwithin{corollary}{section}
%\thispagestyle{empty}
\newcommand{\pr}[1]{\parallel}
\newcommand{\p}[1]{\partial}
%
%

%\begin{frontmatter}
\title{Error and stability analysis of an anisotropic phase-field model for binary-fluid mixtures in the presence of magnetic-field}
\author{A.~Rasheed\thanks{Department of Mathematics, Lahore University of Management Sciences, Opposite Sector U, DHA,Lahore Cantt 54792, Pakistan, amer.rasheed@lums.edu.pk,}, A.~Belmiloudi\thanks{Mathematics Research Institute of Rennes (IRMAR),  European University of Brittany, INSA, 20 Av. Buttes de Coesmes, CS70839, 35708 Rennes Cedex7, France,
aziz.belmiloudi@math.cnrs.fr}}%
\date{}
%
%
%
%
%_____________________________________________________________________________________________________________
%
%                                                   Abstract
%_____________________________________________________________________________________________________________
%
%
%
%
%
%
\maketitle
\begin{abstract}
In this article, we study the error and stability of the proposed numerical scheme in order to solve a two dimensional anisotropic phase-field model with convection and externally applied magnetic field in an isothermal solidification of binary alloys. The proposed numerical scheme is based on mixed finite element method satisfying the CFL condition. A particular application with real physical parameters of Nickel-Copper(Ni-Cu) is considered in order to validate the numerical scheme employed. The results of stability and error analysis substantiates complete accordance with the postulated theoretical error estimates which demonstrates the efficiency of the presented method.

{\bf{keyword}} Phase field methods; Stability and error analysis; binary alloys; magnetic-field.
\end{abstract}
\section{Introduction}
The quality of the final solidified metal is dependent on the dendrite structure evolution during the process of solidification of metals. In order to ameliorate the quality of the solidified materials, it is important to understand and render control over dendritic structure formation during the process of solidification of alloys. During the past decade, the scientists have investigated a great deal of experimental and theoretical studies to explore the microstructure of dendrites during the process of solidification of alloys. Phase field method has been widely used by researchers in order to simulate and study the structure and formation of dendrites (see for instance \cite{Zhu1,Zhu2,AndersonMcFaddenWheeler, BelmiloudiControl,  GrujicicCaoMiller, Ramizer1, Ramizer2, Rappaz, TFT, TonhardtAmberg, TongBeckermann, WarrenBoett}) and the references therein. Unfortunately, the presented models donot render control on dendritic evolution and micro-segregation. Moreover, some experimental observations show that the control on dendrite growth can be achieved in the during the solidification process by applying electric and magnetic fields externally, see for example \cite{LMingjun1, Prescott} and the references therein.
For similar other applications wherin authors have discussed the effect of magnetic field on the metals and alloys, refer to the studies, e.g., for the semi-conductor flow in the melt crystal evolution, \cite{BelmiloudiBook}, for the MHD flows \cite{Hadid1}, \cite{Gunzberger}, \cite{Watanabe}, \cite{Galindo} and for the processes of dendritic solidification,  \cite{Roplekar}, \cite{Prescott},  \cite{Sampath} and the references therein.

In view of aforementioned facts, Rasheed and Belmiloudi  in \cite{RasheedBelmiloudi3} (see also \cite{RasheedBelmiloudi1,RasheedBelmiloudi2, RasheedBelmiloudi4}) has presented a phase-field model which incorporates melt flow and magnetic field.  Initially, the model of Warren-Boettinger \cite{WarrenBoett} is considered and then, among other, Navier Stokes equations, the phase-field and solute equations with a magnetic field applied externally on the domain. The newly developed phase-field model consists three systems, the magnetohydrodynamic system which describes the melt convection by using incompressible Navier-Stokes equations together with the Lorentz force and boussinesq approximations, the phase-field system which represents the phase change and the concentration system describing the concentration change in dendrites during solidification process. The phase field and concentration equations are of convection diffusion type systems. We refer the reader to \cite{RasheedBelmiloudi3} for more detailed  description of the model. 

The existence and uniqueness of the solutions of the derived model is presented in the article \cite{RasheedBelmiloudi3}. In order to perform realistic physical simulations it is indispensable to develop a stable and convergent numerical scheme. This article is devoted to present a numerical approximation scheme using mixed finite-elements and numerical stability and error analysis for Rasheed and Belmiloudi model \cite{RasheedBelmiloudi} for the anisotropic case.  Some numerical results in the isotropic case are presented in the article \cite{AmerNote}.

The paper has been organized as follows. In  section \ref{MF} we recall the mathematical model briefly and its weak formulation in section \ref{WF}. Section \ref{DWF} is describes the discrete variational formulation, in the context of a mixed finite element  method. The numerical stability and error analysis are presented and validated by numerical experiments in section \ref{SANS}.
\section{Mathematical formulation}\label{MF}
Initially a region $\Omega$ is assumed to be occupied by a binary alloy containing two pure elements, the solute B (e.g., Cu) and the solvent A (e.g., Ni), which is electrically conducting incompressible fluid. To construct a numerical scheme and study its numerical stability and error, we recall the model for dendritic solidification given by \cite{RasheedBelmiloudi3}. 
 Let $\mathbf{u}$, $p$, $\psi$, $c$ and $\mathbf{B}$ represent the velocity vector, pressure, phase-field variable, relative concentration and externally magnetic field respectively. Wwe have the following system
\begin{equation}\label{MP}
\hspace{-0.3cm}
\begin{cases}
\displaystyle \rho_{0}\big(\partial_t\mathbf{u} +(\mathbf{u}\cdot\nabla) \mathbf{u}\big)
= 
-\nabla p + \mu \Delta\mathbf{u}+\boldsymbol{\mathcal{A}}_{1}(\psi,c)+b(\psi)\big(\left(\mathbf{u}\times \mathbf{B}\big)\times 
\mathbf{B}\right) & \text{on}~~{\cal Q}, 
\\
\displaystyle\rm{div}\mathbf{u} = 0 & \text{on}~~{\cal Q}, 
\\
\displaystyle\partial_t\psi+(\mathbf{u}\cdot\nabla)\psi 
= \rm{div}\left(\mathcal{A}_{g}(\nabla \psi)\nabla \psi\right) - \mathcal{A}_{2}(\psi,c)
& \text{on}~~{\cal Q},
\\
\displaystyle\partial_t c+(\mathbf{u}\cdot\nabla)c 
=\rm{div}\left(D(\psi)\nabla c +\mathcal{A}_{3}(\psi,c)\nabla\psi\right)
& \text{on}~~ {\cal Q},
\\
\text{subject to the initial conditions}\\
\left(\mathbf{u},\psi,c\right)(t=0)= \left(\mathbf{u}_{0},\psi_{0},c_{0}\right) & \text{on}~~ \Omega,  
\\
\text{and the boundary conditions}\\
\mathbf{u}=0,\quad \mathcal{A}_{g}(\nabla \psi)\nabla \psi\cdot \mathbf{n} =0,
\quad 
(D(\psi)\nabla c +\mathcal{A}_{3}(\psi,c)\nabla\psi) \cdot \mathbf{n}=0 &   \text{on}~~\Sigma,
\end{cases}
\end{equation}
where  $\Omega\subset\Re^{2}$ is a Lipschitz and sufficiently smooth open solidification domain with polygonal boundary $\partial\Omega$, $T$ is the final time of the solidification process, ${\cal Q}=\Omega\times\left(0,T\right)$, $\Sigma=\partial\Omega\times\left(0,T\right)$, $\rho_{0}=\frac{\rho_{0}^{(A)}+\rho_{0}^{(B)}}{2}$ and $\mu=\frac{\mu^{(A)}+\mu^{(B)}}{2}$  are the average density and viscosity,  $D(\psi)$ is the diffusion coefficient and $\mathbf{n}$ is the unit outward normal to $\partial \Omega$. The anisotropy matrix $\mathcal{A}_{g}$ is defined by\\
$$
\mathcal{A}_{g}(\nabla \psi)=M_\psi
\begin{pmatrix}
 \eta^{2}_{\gamma}(\theta) & -\eta_{\gamma}(\theta)\eta'_{\gamma}(\theta)\\  \eta_{\gamma}(\theta)\eta'_{\gamma}(\theta)    &\eta^{2}_{\gamma}(\theta)
\end{pmatrix}
$$
where $M_\psi>0$ and  $\eta_{\gamma}$ is the anisotropy function defined as \cite{WarrenBoett} 
\begin{equation}\label{PFV2}
\eta_{\gamma}=\epsilon_{0}(1+\gamma \cos (k\theta)),
\end{equation}
$\gamma\geq 0$ is the anisotropic amplitude, the integer $k>1$ corresponds to the number of branching directions, $\epsilon_{0}$ is a constant and $\theta$  (the angle between the $x$-axis and $\nabla \psi$)  
\begin{equation}\label{PFV3}
\theta=\arctan\left(\frac{\psi_{y}}{\psi_{x}}\right),
\end{equation}
where $x$ and $y$ are subscripts used to represent the partial derivatives with respect to spatial coordinates, that is $\psi_{x}=\partial \psi/\partial x$ and $\psi_{y}=\partial \psi/\partial y$. For low to moderate accuracy, $\gamma>0$ is selected so that the condition $\displaystyle \eta_{\gamma}(\theta)+\frac{d^{2}\eta_{\gamma}}{d\theta^{2}}(\theta)>0$ is valid for all $\theta$. For $k=4$ (fourfold anisotropy) which is a case of physical appearance, the previous condition is valid if $\gamma<\frac{1}{15}\approx 0,6667$.  
The operators $\boldsymbol{\mathcal{A}_{1}}$, $\mathcal{A}_{2}$ and  $\mathcal{A}_{3}$ are defined by 
\begin{equation}\label{operators}
%\begin{align}
\begin{array}{lcr}
\boldsymbol{\mathcal{A}_{1}}(\psi,c)= \beta_c a_1(\psi) c \mathbf{G}+\zeta {\bf f}(\psi), \quad \displaystyle \mathcal{A}_{2}(\psi,c)
=M_\psi(\frac{\lambda_{1}(c)}{\delta^{2}}g'(\psi)+\frac{\lambda_{2}(c)}{\delta}\bar{p}'(\psi)),\\ 
\displaystyle \mathcal{A}_{3}(\psi,c)
=\alpha_{0}D(\psi)c(1-c)\left(\frac{\lambda_{1}'(c)}{\delta}g'(\psi)-\lambda_{2}'(c)\bar{p}'(\psi)\right), \\
\end{array}
\end{equation}
For a detailed description and derivation of the mathematical model, the reader is referred to \cite{RasheedBelmiloudi3}.
\section{Weak formulation}\label{WF}
The inner product and norm in $L^2(\Omega)$ are defined respectively by
$(.,.)$ and $\mid .\mid$ and 
\begin{align}
\begin{array}{lll}
&\displaystyle{\cal{W}}=(H_{0}^{1}(\Omega))^{2}=\left\{\textbf{v}\in\ \left(H^{1}(\Omega)\right)^{2}\ |\ \textbf{v}=0\ \mbox{on}\ \partial\Omega\right\}, 
\ \ \ \displaystyle{\cal{W}}_{d}=\left\{\textbf{v}\in\ {\cal{W}}\ |\ \rm{div}(\textbf{v})=0\ \right\}, \\  
&\displaystyle{\cal{M}}=H^{1}(\Omega), \ \ \ \displaystyle{\cal{H}}= \left\{q\in L^{2}(\Omega)\ \left| \right. \int_{\Omega} q d\textbf{x} = 0\right\},
\end{array}
\end{align}
where ${\cal{W}}$ is equipped with the norm $\pr\ \nabla .\pr\ $. The scalar product and norm in ${\cal{H}}$ are denoted by the usual $L^2(\Omega)$ inner product and its norm $(.,.)$ and $\mid .\mid$, respectively. 
\begin{Remark}: In order to assure that the pressure is unique, we impose the condition $\displaystyle \int_{\Omega} q d\textbf{x} = 0$ on the pressure which is defined within a class of equivalence, regardless of a time-dependent function. We may impose also other conditions on the pressure, in accordance on its regularity, e.g., the pressure is constant on part of the boundary, etc.
\end{Remark}
The bilinear and trilinear forms are defined as follows (for all $(\mathbf{u},\mathbf{v},\mathbf{w})\in ({\cal{W}})^{3}$,  $p\in {\cal{H}}$, $(c,z)\in ({\cal{M}})^{2}$ and 
$(\phi,\psi)\in ({\cal{M}})^{2}$): 
\begin{align*}
\begin{array}{lll}
&\displaystyle a_u(\mathbf{u},\mathbf{v})
=\mu\int_\Omega\nabla\mathbf{u}\cdot\nabla\mathbf{v} d\mathbf{x}, \ \ c_{p}\left(\textbf{u},p\right)=-\left(\rm{div}(\textbf{u}),p\right), \ \ b_u(\mathbf{u},\mathbf{v},\mathbf{w})= \rho_0 \sum_{i,j=1}^2\int_\Omega u_i(\partial_iv_j)w_j d\mathbf{x},
\\
& \displaystyle b_c(\mathbf{u},c,z)=\sum_{i=1}^2\int_\Omega u_i(\partial_i c) zd\mathbf{x},\ \ b_\psi(\mathbf{u},\psi,\phi)=\sum_{i=1}^2\int_\Omega u_i(\partial_i \psi) \phi d\mathbf{x}.
\end{array}
\end{align*}
Moreover, if $\rm{div}(\mathbf{u})=0$, the trilinear forms satisfy the classical relations given in the following Lemma (see e.g. \cite{BelmiloudiRobin, Temam}):
\begin{Lemma}\label{lem1a}
The trilinear forms $b_{u}, b_{\psi}, b_{c}$ have the following properties
\begin{itemize}
\item [(i)] For all $\textbf{u} \in {\cal{W}}_{d},\ \textbf{v} \in {\cal{W}}$ and $\psi,c \in {\cal{M}}$
\begin{eqnarray}
\nonumber
b_{u}(\textbf{u},\textbf{v},\textbf{v}) = 0,\ \ b_{\psi}(\textbf{u},\psi,\psi) = 0,\ \  b_{c}(\textbf{u},c,c) = 0.
\end{eqnarray}
\item [(ii)] For all $\textbf{u} \in {\cal{W}}_{d},\ \textbf{v}, \textbf{w} \in {\cal{W}}$ and  $\psi,\phi, c,z \in  {\cal{M}}$
\begin{eqnarray}
\nonumber
b_{u}(\textbf{u},\textbf{v},\textbf{w})=-b_{u}(\textbf{u},\textbf{w},\textbf{v}),\ \  
b_{\psi}(\textbf{u},\psi,\phi)=-b_{\psi}(\textbf{u},\phi,\psi),\ \ 
b_{c}(\textbf{u},c,z)=-b_{c}(\textbf{u},z,c).
\end{eqnarray}
\end{itemize}
\end{Lemma}

In order to obtain the weak formulations of the underlying model we multiply the first and second equations of (\ref{MP}) by $\textbf{v} \in {\cal{W}}$, third equation of (\ref{MP}) by $\phi\in {\cal{M}}$ and last equation of (\ref{MP}) by $z\in {\cal{M}}$. Further we integrate the results over $\Omega$ and then using  Green's formulas and boundary conditions, we obtain the following weak formulation of the problem (\ref{MP}) (wherein we added an  artificial source term $\textbf{F}_{u}$, $F_{\psi}$ and $F_{c}$ in each equation of the model for fabricating exact solutions to perform the convergence and stability of the numerical scheme):
Find $(\textbf{u},p,\psi,c) \in \cal{W}\times \cal{H}\times  \displaystyle{\cal{M}}\times  \displaystyle{\cal{M}}$ such that $\forall\ (\textbf{v},q,\varphi,z)\in \cal{W}\times \cal{H}\times  \displaystyle{\cal{M}}\times  \displaystyle{\cal{M}}$
\begin{equation}\label{Ex1MP}
\hspace{-0.3cm}
\begin{cases}
&\displaystyle \rho_{0} \left(\partial_t\textbf{u}, \textbf{v}\right)+a_{u}(\textbf{u},\textbf{v})+b_{u}(\textbf{u},\textbf{u},\textbf{v}) +c_{p}\left(\textbf{v},p\right)- \left(\boldsymbol{\mathcal{A}}_{1}(\psi,c),\textbf{v}\right)
\\
&\displaystyle \hspace{3cm} - \left(b(\psi)((\textbf{u}\times \textbf{B})\times \textbf{B}),\textbf{v}\right) = \left(\textbf{F}_{u},\textbf{v}\right),
\\ 
&\displaystyle -c_{p}\left(\textbf{u},q\right)=0,
\\
&\displaystyle \left(\partial_t\psi,\varphi\right)+b_{\psi}(\textbf{u},\psi,\varphi) +(\mathcal{A}_{g}(\nabla \psi)\nabla\psi,\nabla\varphi) + \left(\mathcal{A}_{2}(\psi,c),\varphi\right)
=\left(F_{\psi},\varphi\right),
\\
&\displaystyle \left(\partial_tc,z\right)+b_{c}(\textbf{u},c,z)+ \left( D(\psi)\nabla c,\nabla z\right)\displaystyle  +  \left(\mathcal{A}_{3}(\psi,c)\nabla \psi,\nabla z\right)=\left(F_{c},z\right),
\\
&\displaystyle \left(\textbf{u},\psi,c\right)(t=0)=\left(\textbf{u}_{0},\psi_{0},c_{0}\right).
\end{cases}
\end{equation}
\section{Discrete weak formulation and finite element discretization}\label{DWF}
Let ${\cal T}_{h}$ be a family of triangulations which discretize  the domain $\overline{\Omega}$ with maximum mesh spacing parameter $0<h=\max_{\mho \in {\cal T}_{h}}diam(\Omega)<h_{0}<1$. To develop the Galerkin approximation of (\ref{Ex1MP}), we consider the $\mathbb{P}_{l}$, $\mathbb{P}_{l-1}$ and $\mathbb{P}_{l}$ finite element subspaces ${\cal W}_{h}$, ${\cal H}_{h}$ and ${\cal M}_{h}$ of $\cal{W}$, $\cal{H}$ and $\cal{M}$ respectively over the partition ${\cal T}_{h}$, where the polynomials $\mathbb{P}_{l}$  is the space of polynomials of degree at most $l$. Furthermore, in order to derive theoretical error estimates of the PDEs like \ref{MP}, we impose the following assumptions (see  \cite{BelmiloudiErrEst, Suli})

\begin{enumerate}

\item[(C1)]  $\exists\ c_{1}>0$, such that $\forall  \textbf{X}=(\textbf{u},\psi,c)\in \left(H^{r+1}(\Omega)\right)^{4}\cap \left({\cal W} \times {\cal M}^{2}\right)$ and $\forall r \in \left[1,l\right]$
\begin{equation*}
%\nonumber
\displaystyle \inf_{\textbf{X}_{h}\in {\cal W}_{h} \times {\cal M}^{2}_{h}} \left\|\textbf{X}-\textbf{X}_{h}\right\|\ \leq\ c_{1} h^{r} \left\|\textbf{X}\right\|_{H^{r+1}(\Omega)}.
\end{equation*}
\item [(C2)] $\exists\ c_{2}>0$, such that  $\forall \ q\in H^{r}(\Omega)\cap {\cal H}$ and $\forall\ r \in \left[1,l\right]$
\begin{equation*}
%\nonumber
\displaystyle \inf_{q_{h}\in {\cal H}_{h}} \left\|q-q_{h}\right\|\ \leq\ c_{2} h^{r} \left\|q\right\|_{H^{r}(\Omega)}.
\end{equation*}
\item [(C3)] $\exists\ c_{3}>0$ such that  (Inf-Sup condition)
\begin{equation*}\label{InfSup}
\displaystyle \inf_{q_{h}\in {\cal H}_{h}} \sup_{\textbf{v}_{h}\in {\cal W}_{h}} \frac{c_{p}(\textbf{v}_{h},q_{h})}{\left\|\textbf{v}_{h}\right\|\left|q_{h}\right|}\ \geq\ c_{3}.
\end{equation*}
\item [(C4)] Let  $\textbf{X}_{0}=(\textbf{u}_{0},\psi_{0},c_{0})\in \left(H^{r+1}(\Omega)\right)^{4}$ with $r \in [1,l]$, then
\begin{equation*}
%\nonumber
\displaystyle  h\left\|\textbf{X}_{0}-\textbf{X}_{0h}\right\|+ \left|\textbf{X}_{0}-\textbf{X}_{0h}\right|\ \leq\ c_{4} h^{r+1},
\end{equation*}
where $\textbf{X}_{0h}=(\textbf{u}_{0h},\psi_{0h},c_{0h})\in {\cal W}_{h}\times {\cal M}^{2}_{h}$ is the approximation of $\textbf{X}_{0}$.
\item [(C5)] For all integers $\mathfrak{m},\mathfrak{p},\mathfrak{q}$ and $\mathfrak{k}$ with $0<\mathfrak{p},\mathfrak{q}\leq \infty$ and $\forall\ K\in {\cal T}_{h}$, we have
\begin{eqnarray*}
%\nonumber
\left\|\textbf{X}_{h}\right\|_{W^{\mathfrak{m},\mathfrak{q}}(K)}\ \leq\ c_{4}\ h^{n/\mathfrak{q}-n/\mathfrak{p}+\mathfrak{k}-\mathfrak{m}}\ \left\|\textbf{X}_{h}\right\|_{W^{\mathfrak{k},\mathfrak{p}}(K)},\ \ \ \forall\ \textbf{X}_{h}\in {\cal W}_{h}\times {\cal M}^{2}_{h}, \\ \nonumber
\left\|\textbf{X}_{h}\right\|_{W^{\mathfrak{m},\mathfrak{q}}(\Omega)}\ \leq\ c_{4}\ h^{n/\mathfrak{q}-n/\mathfrak{p}+\mathfrak{k}-\mathfrak{m}}\ \left\|\textbf{X}_{h}\right\|_{W^{\mathfrak{k},\mathfrak{p}}(\Omega)},\ \ \ \forall\ \textbf{X}_{h}\in {\cal W}_{h}\times {\cal M}^{2}_{h}.
\end{eqnarray*}
\end{enumerate}
The space and time discretization of the problem (\ref{Ex1MP}) i.e. the discrete weak formulation of the problem (\ref{MP}) can now be easily defined. We shall describe in detail the numerical scheme and present the space discretization and the general form of the differential-algebraic systems for (\ref{Ex1MP}). Further we present the time discretization of the problem briefly. The discrete weak formulation is presented as follows: 
Find $(\textbf{u}_{h},p_{h},\psi_{h},c_{h})\in {\cal W}_{h}\times {\cal H}_{h} \times {\cal M}_{h}\times {\cal M}_{h}$ such that  $\forall (\textbf{v}_{h},q_{h},\varphi_{h},z_{h})\in {\cal W}_{h}\times {\cal H}_{h} \times {\cal M}_{h}\times {\cal M}_{h}$ 
\begin{equation}\label{STD}
\hspace{-0.3cm}
\begin{cases}
& \displaystyle\rho_{0}\left(\partial_t\textbf{u}_{h},\textbf{v}_{h}\right)+a_{u}\left(\textbf{u}_{h},\textbf{v}_{h}\right)+b_{u}\left(\textbf{u}_{h},\textbf{u}_{h},\textbf{v}_{h}\right)+ c_{p}\left(\textbf{v}_{h},p_{h}\right) - \left(\boldsymbol{\mathcal{A}}_{1}(\psi_{h},c_{h}),\textbf{v}_{h} \right) \vspace{0.1cm} \\ & \hspace{4cm} - \left( b(\psi_{h})((\textbf{u}_{h}\times \textbf{B})\times \textbf{B}),\textbf{v}_{h}\right)=\left(\textbf{F}_{u},\textbf{v}_{h}\right),\\ 
& \displaystyle-c_{p}\left(\textbf{u}_{h},q_{h}\right)=0, \\ 
& \displaystyle\left(\partial_t\psi_{h},\varphi_{h}\right) + b_{\psi}(\textbf{u}_{h},\psi_{h},\varphi_{h}) +(\mathcal{A}_{g}(\nabla \psi_{h})\nabla\psi_{h},\nabla\varphi_{h}) + \left(\mathcal{A}_{2}(\psi_{h},c_{h}),\varphi_{h}\right)
%& 
=\left(F_{\psi},\varphi_{h}\right),\\
%\ 
& \displaystyle\left(\partial_t c_{h},z_{h}\right)+ b_{c}\left(\textbf{u}_{h},c_{h},z_{h}\right)+\left(D(\psi_{h})\nabla c_{h},\nabla z_{h}\right) +  \left(\mathcal{A}_{3}(\psi_{h},c_{h})\nabla \psi_{h},\nabla z_{h}\right)=\left(F_{c},z_{h}\right),\\
%\
&\left(\textbf{u}_{h},\psi_{h},c_{h}\right)(t=0)=\left(\textbf{u}_{0h},\psi_{0h},c_{0h}\right).
\end{cases}
\end{equation}
Let $(\varphi_{ih})_{1\leq i \leq M}$, $(q_{ih})_{2M+1\leq i \leq 2M+N}$ and $(z_{ih})_{2M+N+1\leq i \leq 2M+N+\tilde{M}}$ be the  basis of ${\cal W}_{h}$, ${\cal H}_{h}$ and ${\cal M}_{h}$ respectively and
\begin{eqnarray}\label{series}\nonumber 
\displaystyle\textbf{u}_{h}(x,t) &=& \sum^{M}_{i=1} \Big[\begin{array}{c}
u_{ih}(t) \\ v_{ih}(t)
\end{array}\Big] \varphi_{ih}(x)=\sum^{M}_{i=1} u_{ih}(t) \underline{\varphi}^{u}_{ih}(x)+\sum^{M}_{i=1} v_{ih}(t) \underline{\varphi}^{v}_{ih}(x),\\  
\displaystyle p_{h}(x,t) &=& \sum^{2M+N}_{i=2M+1} p_{ih}(t) q_{ih}(x),\ \ \ \displaystyle\psi_{h}(x,t) = \sum^{2M+N+\tilde{M}}_{i=2M+N+1} \psi_{ih}(t) z_{ih}(x), \\ \nonumber 
 \displaystyle c_{h}(x,t) &=& \sum^{2M+N+2\tilde{M}}_{i=2M+N+\tilde{M}+1} c_{ih}(t) z_{ih}(x),
\end{eqnarray}
where
$\textbf{u}_{ih} = \left(\begin{array}{c}  u_{ih} \\  v_{ih} \end{array}\right),\ \ \ \underline{\varphi}^{u}_{ih} = \left(\begin{array}{c}  \varphi_{ih} \\  0 \end{array}\right),\ \ \ \underline{\varphi}^{v}_{ih} = \left(\begin{array}{c}  0 \\  \varphi_{ih} \end{array}\right).$\newline 
and $\textbf{v}_{h}$ may be represented as
\begin{eqnarray}\label{testfn}
\textbf{v}_{h} = \Big[\begin{array}{c}
\varphi_{ih}(x) \\ \varphi_{ih}(x) 
\end{array}\Big] =  \Big[\begin{array}{c}
\varphi_{ih}(x) \\ 0
\end{array}\Big] + \Big[\begin{array}{c}
0 \\ \varphi_{ih}(x)
\end{array}\Big]
 = \underline{\varphi}^{u}_{ih}(x) + \underline{\varphi}^{v}_{ih}(x)
\end{eqnarray}
In order to derive the numerical scheme, first few terms have been elaborated subsequently. Making use of (\ref{series}) and (\ref{testfn}), the first term $\rho_{0}\left( \partial_{t}\textbf{u}_{h}, \textbf{v}_{h}  \right)$ in the first equation of weak formulation (\ref{STD}) yields 
\begin{eqnarray}\nonumber 
\rho_{0}\left( \partial_{t}\textbf{u}_{h}, \textbf{v}_{h}  \right) &=& \rho_{0}\left(  \sum^{M}_{i=1} \dfrac{d u_{ih}(t)}{dt} \underline{\varphi}^{u}_{ih}(x)+\sum^{M}_{i=1} \dfrac{d v_{ih}(t)}{dt} \underline{\varphi}^{v}_{ih}(x)\ ,\ \underline{\varphi}^{u}_{ih}(x) + \underline{\varphi}^{v}_{ih}(x) \right) \\ \nonumber  
& = & \rho_{0} \sum^{M}_{i=1} \left( \underline{\varphi}^{u}_{ih}(x)\ ,\ \underline{\varphi}^{u}_{ih}(x)\right) \dfrac{d  u_{ih}(t)}{dt} + \rho_{0}\sum^{M}_{i=1} \left( \underline{\varphi}^{v}_{ih}(x)\ ,\  \underline{\varphi}^{v}_{ih}(x)   \right) \dfrac{d v_{ih}(t)}{dt}
\end{eqnarray}
Consider the second term
\begin{eqnarray}\nonumber 
a_{u}\left( \textbf{u}_{h}, \textbf{v}_{h}  \right) &=& a_{u}\left(  \sum^{M}_{i=1} u_{ih}(t)\underline{\varphi}^{u}_{ih}(x)+\sum^{M}_{i=1} v_{ih}(t) \underline{\varphi}^{v}_{ih}(x)\ ,\ \underline{\varphi}^{u}_{ih}(x) + \underline{\varphi}^{v}_{ih}(x) \right) \\ \nonumber  
& = & \sum^{M}_{i=1} a_{u}\left( \underline{\varphi}^{u}_{ih}(x)\ ,\ \underline{\varphi}^{u}_{ih}(x)\right) u_{ih}(t) + \sum^{M}_{i=1} a_{v}\left( \underline{\varphi}^{v}_{ih}(x)\ ,\  \underline{\varphi}^{v}_{ih}(x)   \right) v_{ih}(t)
\end{eqnarray}
and third term can be obtained in a similar manner as
\begin{eqnarray}\nonumber 
b_{u}\left( \textbf{u}_{h},\textbf{u}_{h},\textbf{v}_{h}   \right) &=& b_{u}\left( \textbf{u}_{h}\  ,\  \sum^{M}_{i=1} u_{ih}(t)\underline{\varphi}^{u}_{ih}(x)+\sum^{M}_{i=1} v_{ih}(t)\underline{\varphi}^{v}_{ih}(x)\ ,\ \underline{\varphi}^{u}_{ih}(x) + \underline{\varphi}^{v}_{ih}(x) \right) \\ \nonumber  
& = & b_{u} \left(\textbf{u}_{h}\ ,\  \sum^{M}_{i=1} u_{ih}(t) \underline{\varphi}^{u}_{ih}(x)\ ,\ \underline{\varphi}^{u}_{ih}(x)\right)  + b_{u}\left(\textbf{u}_{h}\ ,\  \sum^{M}_{i=1}v_{ih}(t) \underline{\varphi}^{v}_{ih}(x)\ ,\  \underline{\varphi}^{v}_{ih}(x)   \right)  \\ \nonumber 
&=& \sum^{M}_{i=1} b_{u}\left(\textbf{u}_{h}\ ,\  \underline{\varphi}^{u}_{ih}(x)\ ,\ \underline{\varphi}^{u}_{ih}(x)\right) u_{ih}(t) + \sum^{M}_{i=1} b_{u}\left(\textbf{u}_{h}\ ,\  \underline{\varphi}^{v}_{ih}(x)\ ,\  \underline{\varphi}^{v}_{ih}(x)   \right) v_{ih}(t)
\end{eqnarray}
Substituting above expressions into the first equation of weak form (\ref{STD}) and using  (\ref{series}) and (\ref{testfn}), the semi-discrete weak form  yields
\begin{eqnarray}
\begin{array}{lll}
\nonumber 
& \displaystyle \sum^{M}_{i=1} \rho_{0}\left(\underline{\varphi}^{u}_{ih},\underline{\varphi}^{u}_{jh}\right)\frac{d u_{ih}}{d t}+\displaystyle \sum^{M}_{i=1}\left(a_{u}\left(\underline{\varphi}^{u}_{ih},\underline{\varphi}^{u}_{jh}\right)+b_{u}\left(\textbf{u}_{h},\underline{\varphi}^{u}_{ih},\underline{\varphi}^{u}_{jh}\right) \right. \vspace{0.2cm} \\ & \hspace{2cm} \left. - \left(b(\psi_{h})((\underline{\varphi}^{u}_{ih}\times \textbf{B})\times \textbf{B}), \underline{\varphi}^{u}_{jh}\right) \right) u_{ih}  + \displaystyle  \sum^{2M+N}_{i=2M+1} \left(q_{ih},div(\underline{\varphi}^{u}_{jh})\right)p_{ih} \\ &  \hspace{2cm}  -  \left(\boldsymbol{\mathcal{A}}_{1}(\psi_{h},c_{h}),\underline{\varphi}^{u}_{jh}\right) + \displaystyle  \sum^{M}_{i=1}\rho_{0}\left(\underline{\varphi}^{v}_{ih},\underline{\varphi}^{v}_{jh}\right)\frac{d v_{ih} }{d t} \\ &  \hspace{2cm}  +\displaystyle \sum^{M}_{i=1}\left(a_{u}\left(\underline{\varphi}^{v}_{ih},\underline{\varphi}^{v}_{jh}\right)+b_{u}\left(\textbf{u}_{h},\underline{\varphi}^{v}_{ih},\underline{\varphi}^{v}_{jh}\right) \right. \vspace{0.2cm} \\ &  \hspace{2cm} \left.  - \left(b(\psi_{h})((\underline{\varphi}^{v}_{ih}\times \textbf{B})\times \textbf{B}), \underline{\varphi}^{v}_{jh}\right) \right) v_{ih}+ \displaystyle  \sum^{2M+N}_{i=2M+1} \left(q_{ih},div(\underline{\varphi}^{v}_{jh})\right)p_{ih}  \vspace{0.1cm} \\ & \hspace{2cm}  - \left(\boldsymbol{\mathcal{A}}_{1}(\psi_{h},c_{h}),\underline{\varphi}^{v}_{jh}\right)
  =\left(\textbf{F}_{u},\underline{\varphi}^{u}_{jh} + \underline{\varphi}^{v}_{jh} \right),\ \ \ \text{~for all~} 1\leq j \leq M,\ \ \ \  \vspace{0.5cm}  \\ 

& \displaystyle -\sum^{2M+N}_{i=2M+1} \Bigl\{ \left(div(\underline{\varphi}^{u}_{ih}),q_{jh}\right)u_{ih} +\left(div(\underline{\varphi}^{v}_{ih}),q_{jh}\right)v_{ih} \Bigr\} =0,\ \ \ 2M+1\leq j \leq 2M+N,  \vspace{0.5cm}  \\ 

& \displaystyle \sum^{2M+N+\tilde{M}}_{i=2M+N+1}\left(z_{ih},z_{jh}\right)\frac{d \psi_{ih}}{d t}+\sum^{2M+N+\tilde{M}}_{i=2M+N+1}\Bigl\{ b_{\psi}\left(\textbf{u}_{h},z_{ih},z_{jh}\right) 
+\left(\mathcal{A}_{g}(\nabla \psi_{h})\nabla z_{ih},\nabla z_{jh}\right) \Bigr\}  \psi_{ih}   \vspace{0.2cm}  \\ & \hspace{2cm}  + \left(A_{2}\left(\psi_{h},c_{h}\right),z_{jh}\right)  \displaystyle   =\left(F_{\psi},z_{jh}\right),\ \ \text{~for all~} 2M+N+1\leq j \leq 2M+N+\tilde{M}, \vspace{0.5cm}  \\ 

& \displaystyle \sum^{2M+N+2\tilde{M}}_{i=2M+N+\tilde{M}+1}\left(z_{ih},z_{jh}\right)\frac{d c_{ih}}{d t}+ \sum^{2M+N+2\tilde{M}}_{i=2M+N+\tilde{M}+1}\Bigl\{ b_{c}\left(\textbf{u}_{h},z_{ih},z_{jh}\right)  \vspace{0.2cm}  \\ &  \hspace{2cm}  + \left(D(\psi_{h})\nabla z_{ih},\nabla z_{jh}\right)\Bigr\} c_{ih} \displaystyle  + \sum^{2M+N+\tilde{M}}_{i=2M+N+1} \left(A_{3}(\psi_{h},c_{h})\nabla z_{ih},\nabla z_{jh}\right)\psi_{ih}  \vspace{0.1cm}   \\ & \hspace{2cm}  \displaystyle   =\left(F_{c},z_{jh}\right),\ \ \text{~for all~} 2M+N+\tilde{M}+1\leq j \leq 2M+N+2\tilde{M}.
\end{array}
\end{eqnarray}
These equations can be written in the differential-algebraic system form (DAE) as 
\begin{align}\label{NumScheme0}
&\mathbb{M}\ \frac{d \mathbf{Y}^{h}}{d t} + \mathbb{A}(\mathbf{Y}^{h})\mathbf{Y}^{h} + \mathbf{L}(\mathbf{Y}^{h})= \mathbf{R}, 
\qquad 
\mathbf{Y}^{h}(t=0)=\mathbf{Y}^{h}_{0},
\\ \nonumber 
&\mathbf{Y}^{h}=
\begin{pmatrix}
\mathbf{u}_{1h}\,\cdots\,\mathbf{u}_{Mh}&p_{1h}\,\cdots\,p_{Nh}&
\psi_{1h}\,\cdots\,\psi_{\tilde{M}h}&c_{1h}\,\cdots \,c_{\tilde{M}h}
\end{pmatrix}^{t},
\end{align}
where 
$\mathbf{R}= \begin{pmatrix} R_{1} & 0 & R_{3} & R_{4} 
\end{pmatrix}^{t}$,  $\mathbf{L}(Y_{h})= \begin{pmatrix}
L_{1} & 0 & L_{3} & 0 \end{pmatrix}^{t}$ and, for $ K_1=2M+N+2\tilde{M}$ and $K_2=2M+N+2\tilde{M}$
\begin{align*}
&\mathbb{M}=
\begin{pmatrix}
M_{11} & 0 & 0 & 0 \\  0 & 0 & 0 & 0 \\  0 & 0 & M_{33} & 0 \\  0 & 0 & 0 & M_{44}
\end{pmatrix}\in\mathbb{R}^{K_1,K_2},
\qquad
\mathbb{A}(\mathbf{Y}_{h})=
\begin{pmatrix} 
A_{11} & A_{12} & 0 & 0 \\  A_{21} & 0 & 0 & 0 \\  0 & 0 & A_{33} & 0 \\  0 & 0 & A_{43} & A_{44} 
\end{pmatrix}
 \in\ \mathbb{R}^{K_1,K_2},
\end{align*}
with
\begin{align}\label{Matrix}
\begin{array}{lll}
&\displaystyle \left(M_{11}\right)_{ji}=\rho_{0}\left(\underline{\varphi}^{u}_{ih},\underline{\varphi}^{u}_{jh}\right)+\rho_{0}\left(\underline{\varphi}^{v}_{ih},\underline{\varphi}^{v}_{jh}\right),
\quad 
\left(M_{33}\right)_{ji} = \left(z_{ih},z_{jh}\right),
\quad 
\left(M_{44}\right)_{ji} = \left(z_{ih},z_{jh}\right),
\\
&\left(A_{11}\right)_{ji}=a_{u}\left(\underline{\varphi}^{u}_{ih},\underline{\varphi}^{u}_{jh}\right) +a_{u}\left(\underline{\varphi}^{v}_{ih},\underline{\varphi}^{v}_{jh}\right) + b_{u}\left(\mathbf{u}_{h},\underline{\varphi}^{u}_{ih},\underline{\varphi}^{u}_{jh}\right) + b_{u}\left(\mathbf{u}_{h},\underline{\varphi}^{v}_{ih},\underline{\varphi}^{v}_{jh}\right)
\\ 
& \hspace{1.5cm} -\left(b(\psi_{h})((\underline{\varphi}^{u}_{ih}\times \mathbf{B})\times \mathbf{B}), \underline{\varphi}^{u}_{jh}\right)-\left(b(\psi_{h})((\underline{\varphi}^{v}_{ih}\times \mathbf{B})\times \mathbf{B}), \underline{\varphi}^{v}_{jh}\right),
\\
& \left(A_{12}\right)_{ji}=\left(q_{ih},\nabla\cdot(\underline{\varphi}^{u}_{jh})\right)
+\left(q_{ih},\nabla\cdot\underline{\varphi}^{v}_{jh}\right)
=\left(A_{21}\right)_{ij},
\quad\,\\
&\left(A_{33}\right)_{ji}=\left(\mathcal{A}_{g}(\nabla \psi_{h})\nabla z_{ih},\nabla z_{jh}\right) + b_{\psi}\left(\mathbf{u}_{h},z_{ih},z_{jh}\right),
\\
&\left(A_{43}\right)_{ji}=\left(\mathcal{A}_2(\psi_{h},c_{h})\nabla z_{ih},\nabla z_{jh}\right),
\qquad
\left(A_{44}\right)_{ji}=\left(D(\psi_{h})\nabla z_{ih},\nabla z_{jh}\right) + b_{c}\left(\mathbf{u}_{h},z_{ih},z_{jh}\right),
\\
&\left(L_{1}\right)_{j}=\left(\boldsymbol{\mathcal{A}}_{1}(\psi_{h},c_{h}),\underline{\varphi}^{u}_{jh}\right) +\left(\boldsymbol{\mathcal{A}}_{1}(\psi_{h},c_{h}),\underline{\varphi}^{v}_{jh}\right), 
\qquad
\left(L_{3}\right)_{j}=\epsilon_1\left(\mathcal{A}_2(\psi_{h},c_{h}),z_{jh}\right),
\\
&\left(R_{1}\right)_{j}=\left(\mathbf{F}_{u},\underline{\varphi}^{u}_{jh}\right) +\left(\mathbf{F}_{u},\underline{\varphi}^{v}_{jh}\right),
\quad
\left(R_{3}\right)_{j}=\left(F_{\psi},z_{jh}\right),
\quad \,\,
\left(R_{4}\right)_{j}=\left(F_{c},z_{jh}\right).
\end{array}
\end{align}
The equation (\ref{NumScheme0}) can be written in general form as
\begin{equation}\label{NumScheme}
\displaystyle {\cal F}(t,\mathbf{Y}^{h}(t),\frac{d \mathbf{Y}^{h}}{dt})= 0,\ \ \ \ \mathbf{Y}^{h}(t=0)=\mathbf{Y}^{h}_{0}.
\end{equation}

In order to consider the fully discrete scheme, for an integer $K>0$, we introduce the timestep  $\tau=\frac{T}{K}$, the time subdivision $t_{i}=i\tau$ ($0 \leq i \leq K$) of $[0,T]$ and, for sufficiently regular function $\textbf{v}$, we denote by $\textbf{v}_{i}$ the value of $\textbf{v}$ at time $t_{i}$ and by $\p\ _{\tau,n}\textbf{v}=\frac{\textbf{v}_{n+1}-\textbf{v}_{n}}{\tau}$.

The differential-algebraic system (\ref{NumScheme0}) is first fully discretized by invoking Euler's backward difference method as 
\begin{equation}\label{FNumScheme}  
\displaystyle {\cal F}(t_{n+1},\mathbf{Y}^{h}_{n+1},\p\ _{\tau,n}\mathbf{Y}^{h})= 0
\end{equation}
and then solved by employing the Newton iteration technique on the resulting non-linear fixed-point system, for this purpose we have used the solver DASSL \cite{LindaPetzold}.

Then, the following a priori error estimates for the solution $(\Psi_{h}, p_{h})$, with $\Psi_{h}=(\textbf{u}_{h},\psi_{h},c_{h})$, of the finite element discretization (\ref{FNumScheme}) can be
obtained by using the assumptions (C1)-(C5) and method as those in e.g. \cite{BelmiloudiErrEst, Suli} (for some $\beta_{1}, \beta_{2}>1$ and $\alpha\geq 1$) 
\begin{equation}\label{CF1}
\left\|\Psi_{h}-\Psi\right\|_{\ell^{2}(0,T,L^{2}(\Omega))} \leq C_{\delta}(\tau^{\alpha}+h^{\beta_{1}}) \text{~and~}
\left\|p_{h}-p\right\|_{\ell^{2}(0,T,L^{2}(\Omega))} \leq C_{\delta}(\tau^{\alpha}+h^{\beta_{2}})
\end{equation}
where $C_{\delta}>0$ is independent of $h$, $\Psi=(\textbf{u},\psi,c)$ is the known exact solution of the problem under consideration and  the space $\ell^{p}(0,T,{\mathbb{X}})$, for a Banach space ${\mathbb{X}}$ and $0<p <+\infty$ is defined by 
\begin{eqnarray}
\nonumber
\ell^{p}(0,T,{\mathbb{X}})= \Bigl\{\textbf{v}:(t_{1},...,t_{k})\rightarrow {\mathbb{X}}\ \mbox{such that}\ \left\|\textbf{v}\right\|_{\ell^{p}(0,T,\mathbb{X})}=\left(\tau \sum^{k}_{i=1}\left\|\textbf{v}_{i}\right\|^{p}_{\mathbb{X}}\right)^{1/p}<\infty\Bigr\}. \\ \nonumber 
\end{eqnarray}
%
%It is to be noted that $\beta_{i}$, $i=1,2,$ are greater than $1$ and less than minimum of the degree of the finite elements (polynomials) and the Sobolev space regularity of the solutions. Moreover, for optimal spatial (resp. temporal) convergence rate we take $\tau^{\alpha}\leq h^{\beta_{i}}$, $i=1,2$ (resp. $h^{\beta_{i}}\leq \tau^{\alpha}$, $i=1,2$).

In the next section, we shall present the results (\ref{CF1}) in order to validate the error analysis and stability of numerical scheme by by considering some examples.
\section{Analysis of the numerical scheme}\label{SANS}
We shall discuss error analysis of the scheme with the help of numerical examples to verify theoretical estimates (\ref{CF1}) and stability of the method. To obtain the convergence rates of the derived scheme, two numerical tests are conducted: the first compute the time discretization error, and the second calculate 
the spatial discretization error. In order to validate the numerical stability of the method, we incorporated $(1 -\epsilon\ randf)$ on the right hand side terms to
introduce perturbations by introducing $\epsilon$ in the numerical solution, where $randf$ denotes random function (that generates some F-distributed random variables) assumes values in $[0, 1]$ and $\epsilon$ is parameter to control perturbation. The parameters are same as given in \cite{WarrenBoett} and the constants for the melt–flow equations are taken in view of the physical properties of the nickel–copper (Ni–Cu) system see Table \ref{PhysicalConstants1} (see e.g., \cite{RasheedBelmiloudi3}). 

\begin{table}[!ht]
%\vspace{-0.5cm}
\small
\begin{center}
\begin{tabular}{ l | c | c | c | c } \hline\hline Property Name &  Symbol & Unit & Nickel (A) &  Copper (B) \\ \hline \hline Melting temperature  &  $T_{m}$ & $K$ & $1728$ & $1358$ \\ \hline Latent heat  & $L$  & $J/m^{3}$ & $2350\times 10^{6}$ & $1758\times 10^{6}$ \\ \hline   Diffusion coeff. liquid & $D_{L}$ & $m^{2}/s$ & $10^{-9}$ & $10^{-9}$ \\ \hline  Diffusion coeff. solid & $D_{S}$ & $m^{2}/s$ & $10^{-13}$ & $10^{-13}$ \\ \hline Linear kinetic coeff.  & $\beta$ & $m/K/s$ & $ 3.3 \times 10^{-3}$ & $3.9 \times 10^{-3}$ \\ \hline Interface thickness & $\delta$ & $m$ & $8.4852 \times 10^{-8}$ & $6.0120 \times 10^{-8}$ \\ \hline Density & $\rho$ & $Kg/m^{3}$ & $7810$ & $8020$ \\ \hline  viscosity & $\mu$ & $Pa\cdot s$ & $4.110\times 10^{-6}$ & $0.597\times 10^{-6}$ \\ \hline Surface energy & $\sigma$ & $J/m^{2}$ & $0.37$ & $0.29$ \\ \hline Electrical conductivity & $\sigma_{e}$ & $S/m$ & $14.3\times 10^{6}$ & $59.6\times 10^{6}$ \\ \hline Molar volume & $V_{m}$ & $m^{3}$ & $7.46\times 10^{-6}$ & $7.46\times 10^{-6}$ \\  \hline Mode Number & $k$ & N/A & 4 & 4  \\ \hline Anisotropy Amplitude & $\gamma_{0}$ & N/A & $0.04$ & $0.04$ \\ \hline \hline \end{tabular} \caption{Physical values of constants} \label{PhysicalConstants1}
\end{center}
\vspace{-0.5cm}
\end{table}

\begin{figure}[htb]
\begin{center}
\includegraphics[width=10cm,height=10cm]{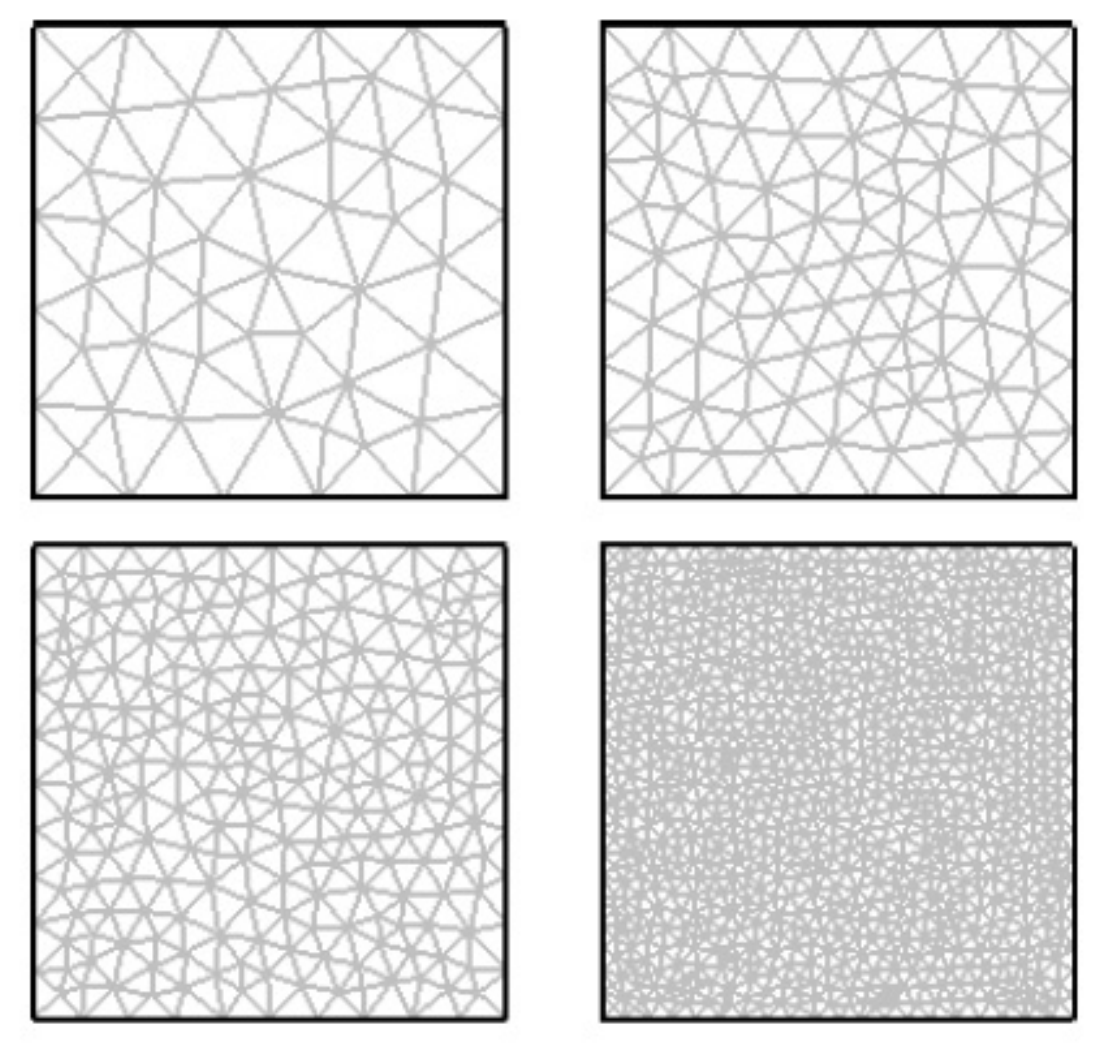} 
\caption{Meshes}
\label{Mesh}
\end{center}
\end{figure}

As exact solutions, we consider the two following examples (with $T=1$ and ${\bf B}=\frac{1}{\sqrt{2}}(1,1)$). 
\begin{itemize}
\item \underline{Example 1} :
\begin{equation}\label{EX1}
\hspace{-0.5cm}
\begin{array}{lll}
\displaystyle u_{ex}(x,y,t)=\frac{2}{(2\pi)^{2}}e^{1-t}sin(x)^{2}y(1-\frac{y}{2\pi})(1-\frac{y}{\pi}),\\ 
\displaystyle v_{ex}(x,y,t)=-\frac{2}{(2\pi)^{2}}e^{1-t}sin(x)cos(x)y^{2}(1-\frac{y}{2\pi})^{2},~~p_{ex}(x,y,t)=  e^{1-t}cos(y),\\ 
\displaystyle \psi_{ex}(x,y,t)=\frac{e^{1-t}}{2}(cos(x)cos(y)+1),~~c_{ex}(x,y,t)=  \frac{8}{(2\pi)^{2}}e^{1-t}x^{^2}(1-\frac{x}{2\pi})^{2}(cos(y)+1), 
\end{array}
\end{equation}
where $\Omega=[0,2\pi]\times[0,2\pi]$. 
\item \underline{Example 2} :
\begin{equation}\label{EX2}
\hspace{-0.5cm}
\begin{array}{lll}
\displaystyle u_{ex}(x,y,t)=  4\pi e^{t-1} x^{2}(1-x)^{2} sin(2 \pi y) cos(2 \pi y),\\
\displaystyle v_{ex}(x,y,t)= -2 e^{t-1}x(2x^{2}-3x+1)sin^{2}(2\pi y),~~p_{ex}(x,y,t)=e^{t-1}cos(2\pi x),\\
\displaystyle \psi_{ex}(x,y,t)=  \frac{1}{4} e^{t-1}(cos(2 \pi x)+cos(2\pi y)+2), ~~c_{ex}(x,y,t)= 8 e^{t-1}(x^{2}(1-x)^{2}+y^{2}(1-y)^{2}),
%\nonumber
\end{array}
\end{equation}
%\end{eqnarray}
where $\Omega=[0,1]\times[0,1]$. 
\end{itemize}
The terms $\textbf{F}_{u}$, $F_{\psi}$ and $F_{c}$ on right hand sides are computed carefully to ensure that (\ref{EX1}) (resp. (\ref{EX2})) is the exact solution of system (\ref{MP}). 
We consider four meshes with step size h (see Fig.\ref{Mesh} and Table \ref{MeshStat}).
\begin{table}[!ht]
\begin{center}
\begin{tabular}{c|c|c|c} \hline\hline Mesh~~~~~~~~~~ &  $h$~~~~~~~~ & Elements~~~~~~~~ & Boundary elements\\ \hline\hline 1.0 & 0.20 & 106.0 & 5.0 \\ \hline 2.0 & 0.150 & 200.0 & 7.0 \\ \hline 3.0 & 0.10 & 434.0 & 10.0 \\ \hline 4.0 & 0.050 & 1712.0 & 20.0 \\ \hline \hline \end{tabular} \caption{Mesh Statistics}
\label{MeshStat} 
%5 & 0.01 & 42904 & 100 \\ \hline
\vspace{-1cm}
\end{center}
\end{table}
\subsection{Numerical error analysis}\label{ErrAna}
In order to validate the error estimates numerically and the convergence orders of scheme, two types of computations have been made. First, we have computed the convergence rates with respect to apatial coordinates wherein sufficiently small timesteps $\tau$ (as compared to the spatial step size $h$) are fixed and we have varied the spatial step size $h$ as described in Table \ref{MeshStat} of mesh statistics. The rates $\beta_{1}$ and $\beta_{2}$ have been calculated with respect to $h$ and we use the Lagrange-quadratic $\mathbb{P}_{2}$ and Lagrange-cubic $\mathbb{P}_{3}$ finite elements for the phase-field and concentration system, and the velocity/pressure mixed finite elements $\mathbb{P}_{2}-\mathbb{P}_{1}$ and $\mathbb{P}_{3}-\mathbb{P}_{2}$  for the melt flow system. 

\begin{figure}[!hb]
\centering
\subfigure[Spatial error curves in Example 1]{
   \includegraphics[width=6.5cm,height=5cm]{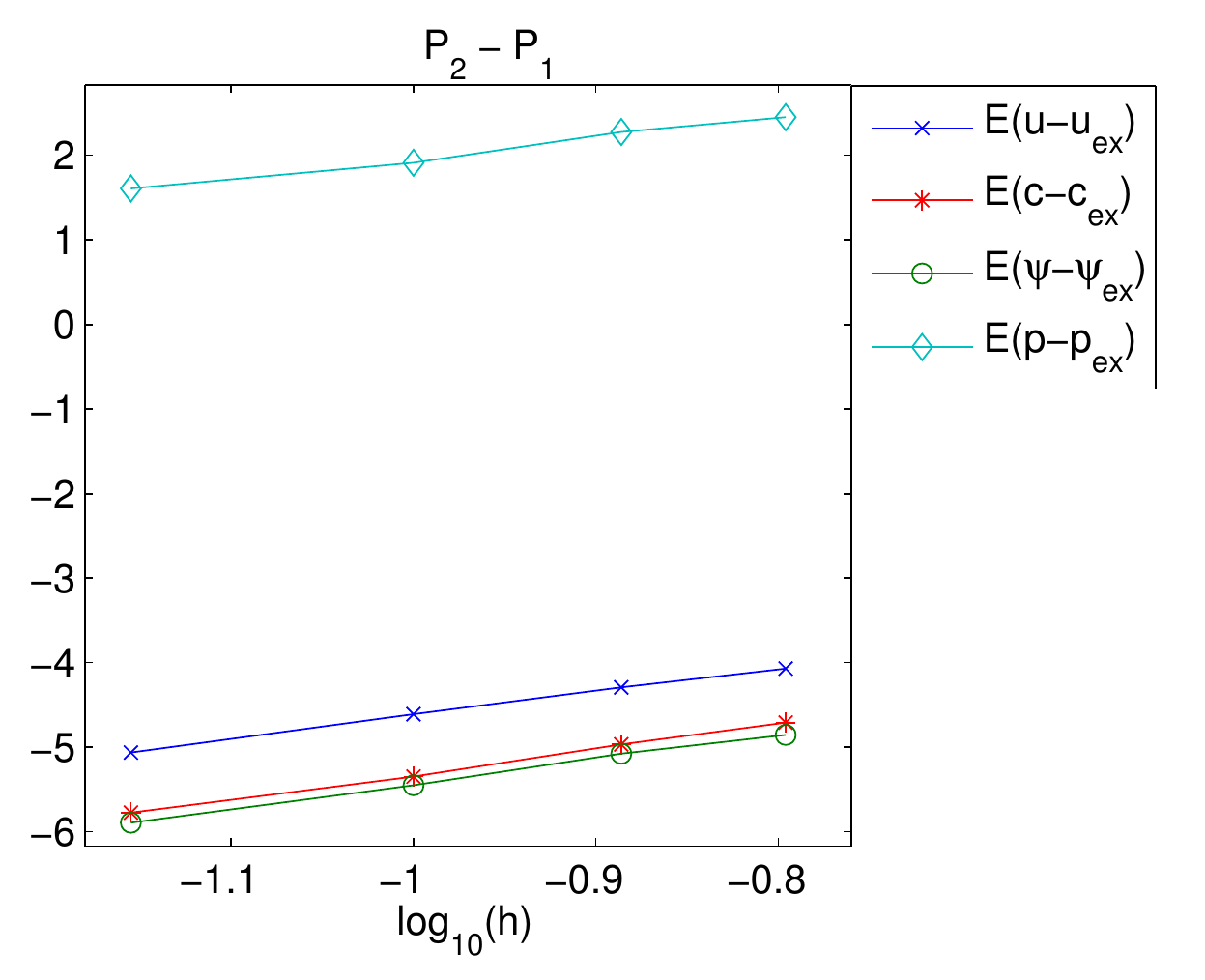}
   \includegraphics[width=6.5cm,height=5cm]{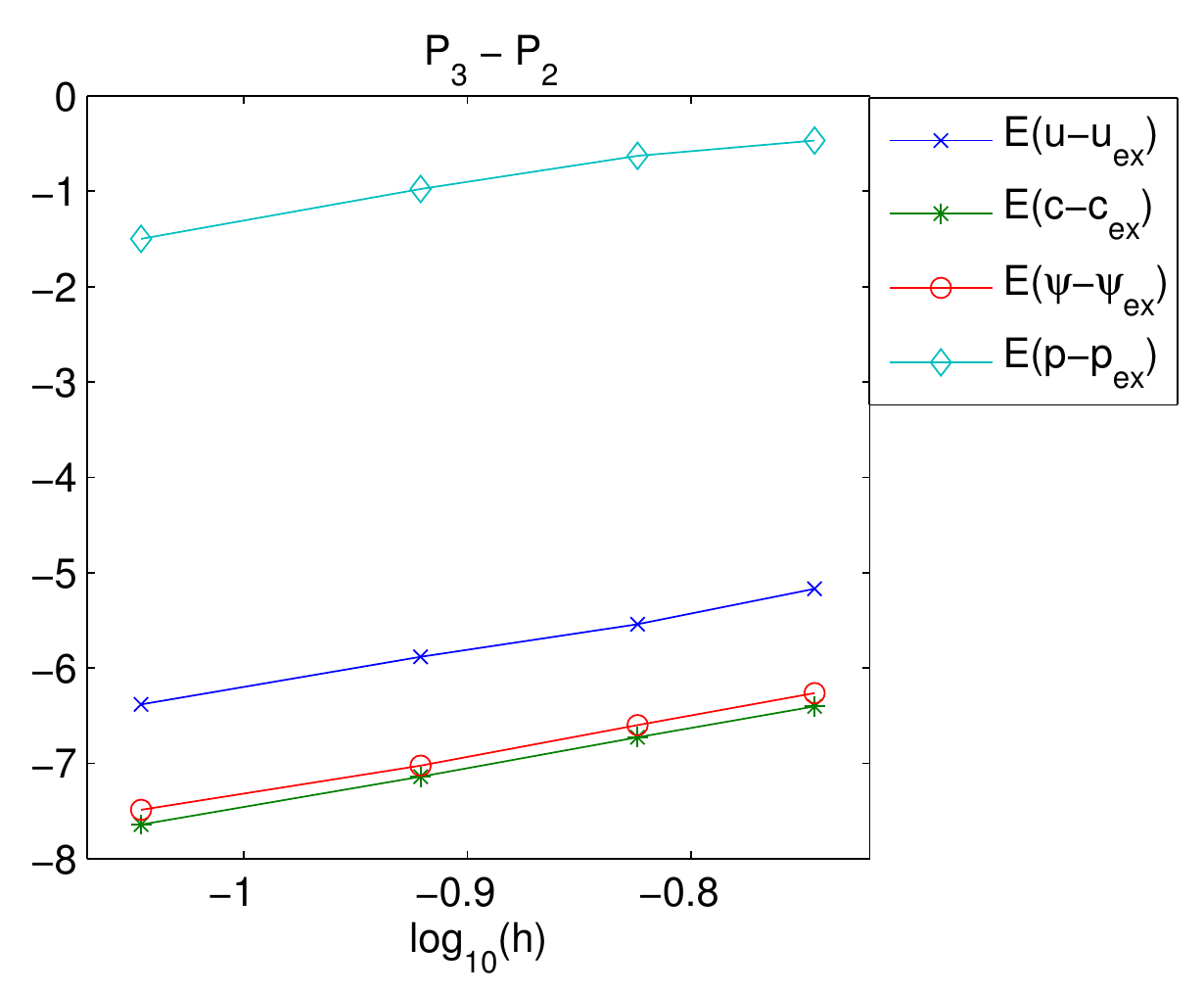}
   \label{Fig:ErrorSpacea}}
 \subfigure[Spatial error curves in Example 2]{ \includegraphics[width=6.5cm,height=5cm]{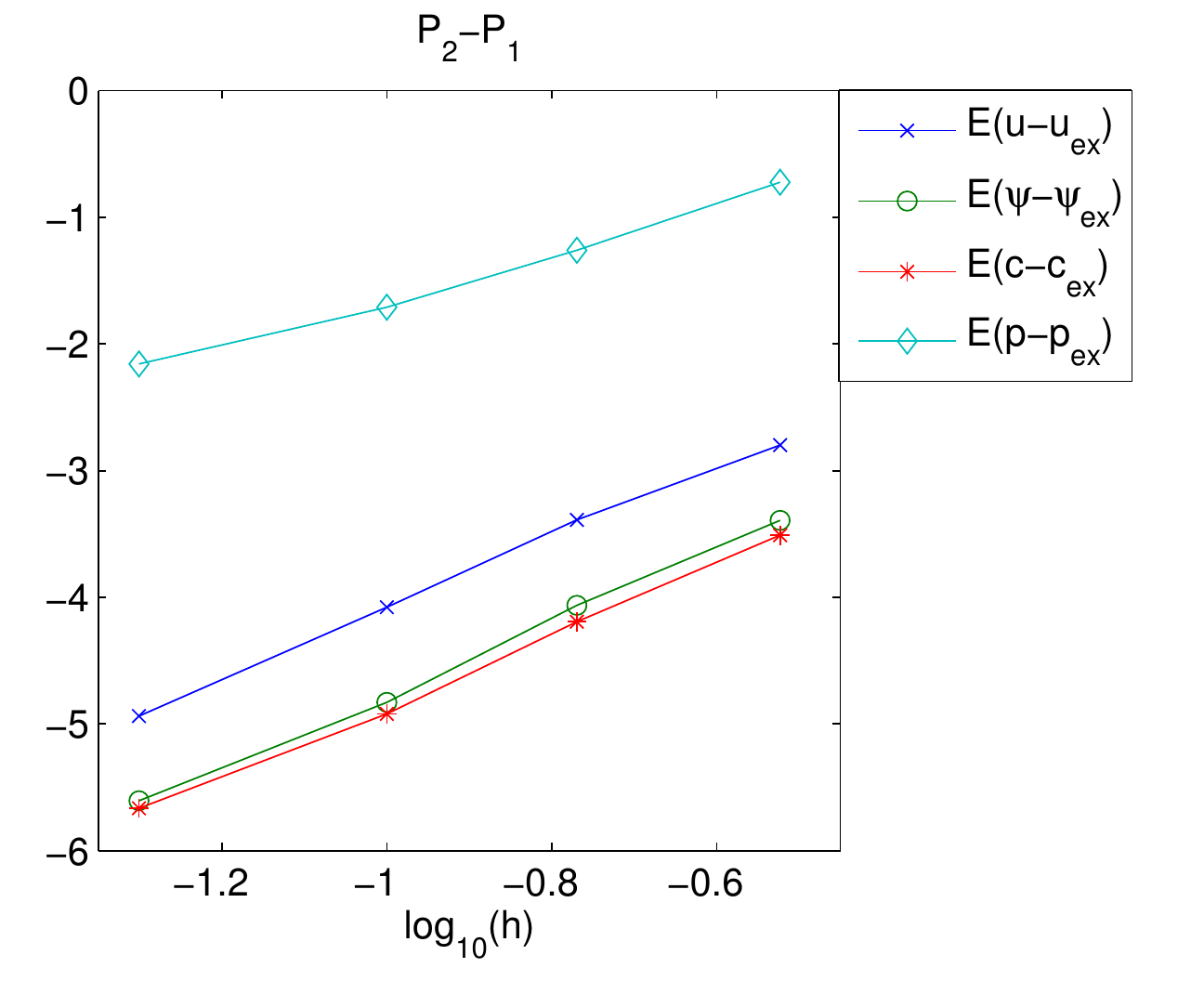}
   \includegraphics[width=6.5cm,height=5cm]{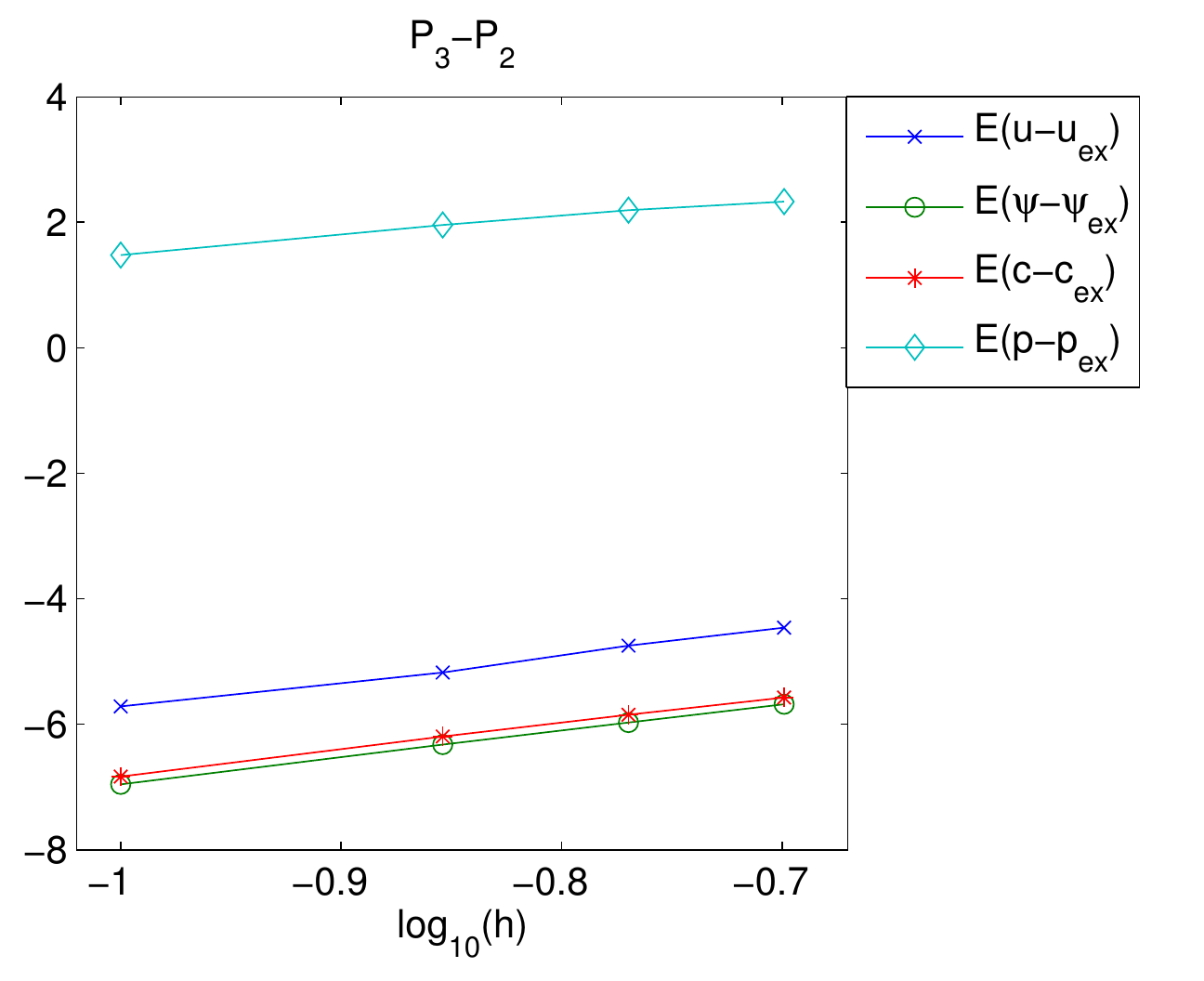}
   \label{Fig:ErrorSpaceb}}
\caption{Error curves with respect to spatial step $h$ obtained in Examples $1$ and $2.$}\label{Fig:ErrorSpace}
\end{figure}

\begin{figure}[!ht]
\centering
\subfigure[Temporal error curves in Example 1]{
   \includegraphics[width=6.5cm,height=5cm]{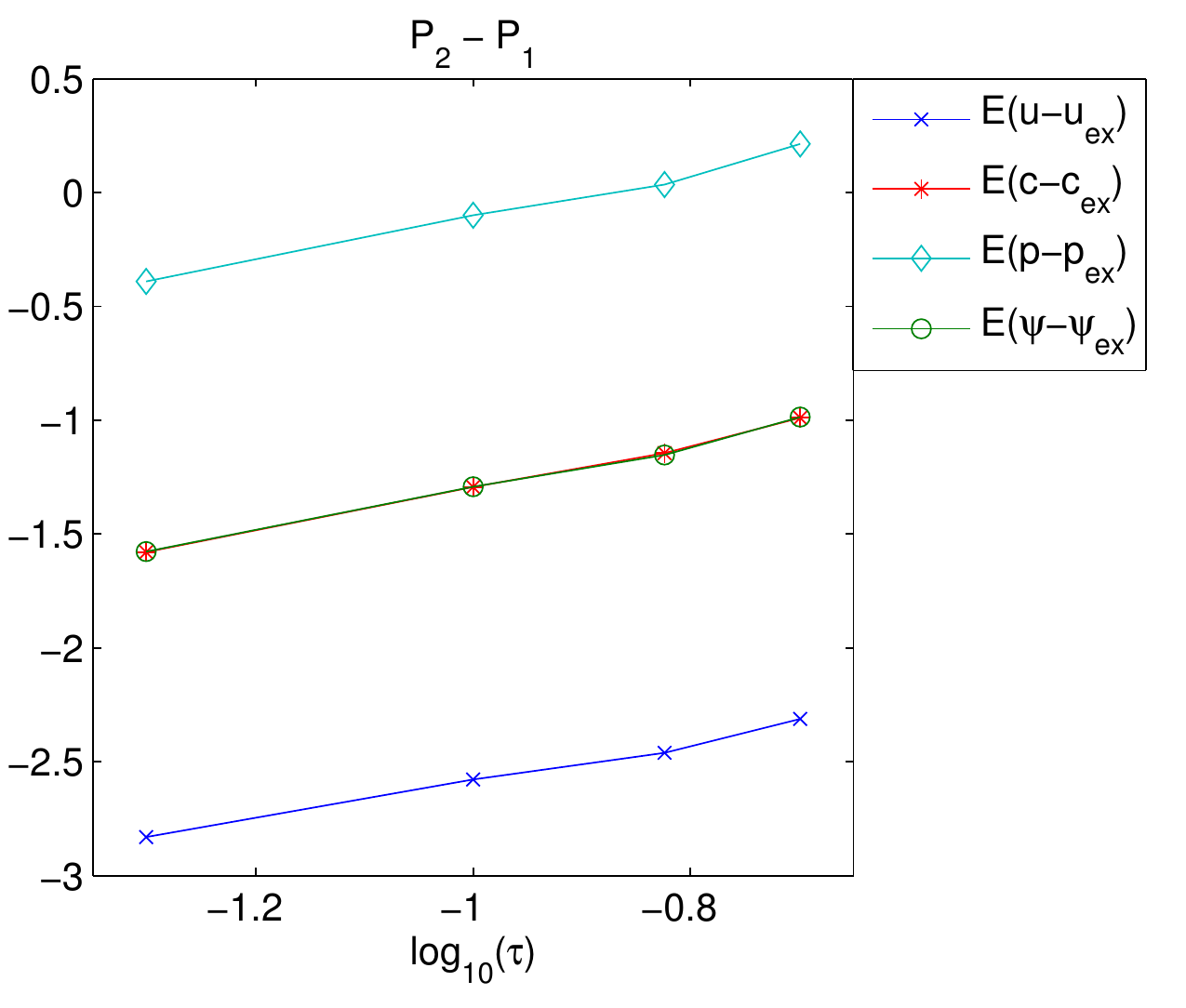}
   \includegraphics[width=6.5cm,height=5cm]{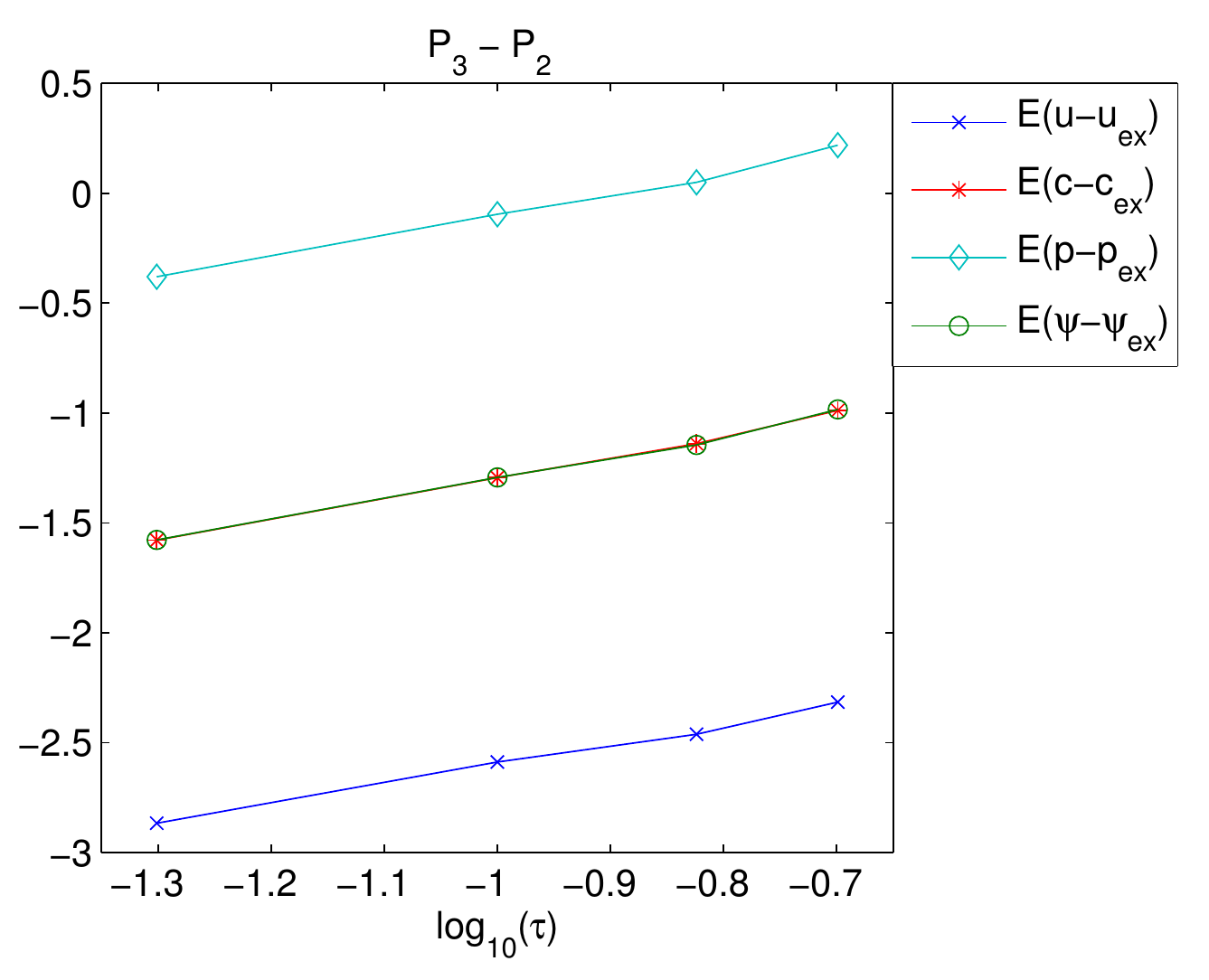}
   \label{fig:Ex2ErrorTimea}}
 \subfigure[Temporal error curves in Example 2]{ \includegraphics[width=6.5cm,height=5cm]{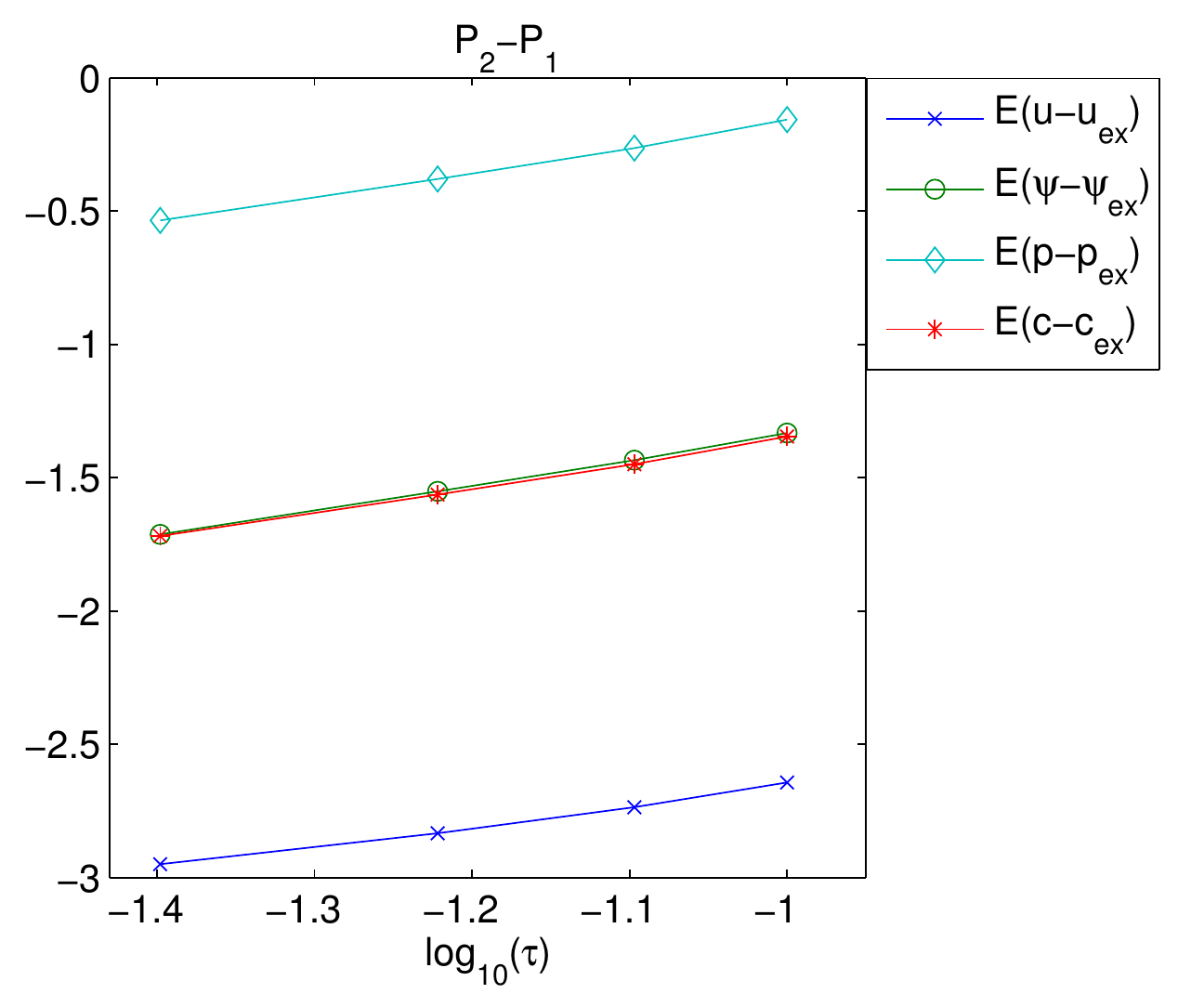}
   \includegraphics[width=6.5cm,height=5cm]{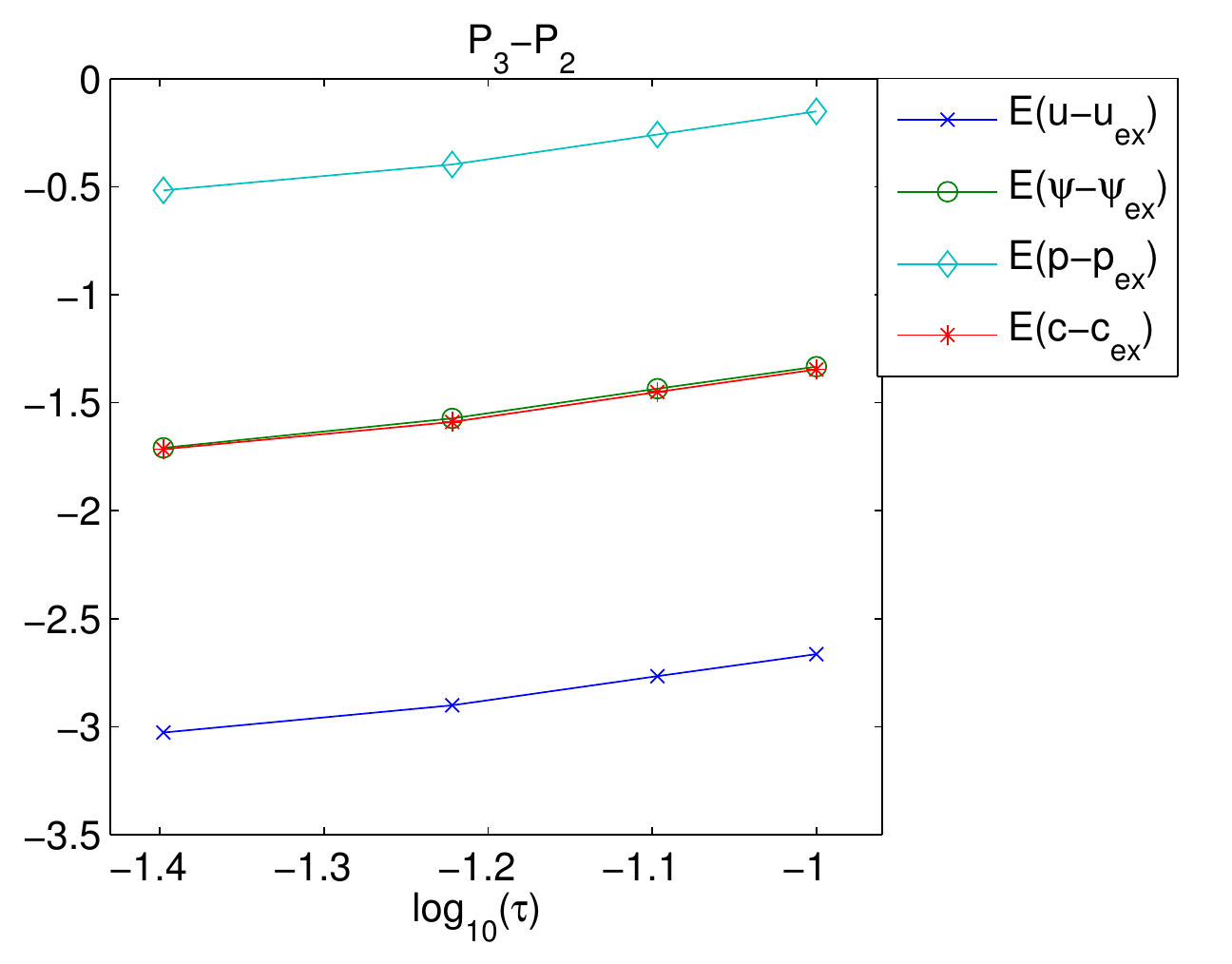}
   \label{fig:Ex2ErrorTimeb}}
\caption{Error curves with respect to temporal step $\tau$ obtained in Examples $1$ and $2$.}\label{ErrorTime}
\end{figure}

The Fig. \ref{Fig:ErrorSpace} and Fig. \ref{ErrorTime} represent the $L^{2}({\cal Q})$-norms of errors  of $\textbf{u}$, $p$, $\psi$ and $c$ which are plotted versus $h$ and $\tau$ respectively, in $log$-scales. For $h$-curves, we use  $\tau=0.01$, $\tau=0.001$ and $\tau=0.0001$ for  linear, quadratic and cubic polynomials, respectively. It is to be noted that the slopes of error curves  are approximately equal to $3$ and $4$ for quadratic and cubic finite elements respectively for the velocity, phase-field and concentration, whereas the slopes of error curves for the pressure are approximately equal to $2$ and $3$ for linear and quadratic finite elements respectively; refer to Table \ref{tab:herrbeta1} and Table \ref{tab:herrbeta2}. $\tau$-curves slopes of all the curves are approximately 1, i.e., $\alpha=1$; refer to 
Table \ref{tab:tauerralpha1} and Table \ref{tab:tauerralpha2}. Both numerical estimates are in an excellent agreement with error estimates (\ref{CF1}).\newline

\begin{table}[!ht]
\centering % centering table
\begin{tabular}{ c | c | c | c }
 \hline\hline % inserting double-line
   & Error Estimate & $\textbf{P}_{2}-\textbf{P}_{1}$ & $\textbf{P}_{3}-\textbf{P}_{2}$
 \\ [0.5ex]
\hline    % inserts single-line
% Entering 1st row
& $\beta_{1}$ for $\textbf{u}$ & $2.6201$ & $3.8730$ \\ [-0.7ex] 
\raisebox{1.5ex}{Example 1} & $\beta_{2}$ for $p$ & $1.9207$ & $3.0646$  \\ [0.5ex]  \hline 
% inserts single-line
% Entering 1st row
&  $\beta_{1}$ for $\textbf{u}$ & $2.7664$ & $4.0303$ \\ [-0.7ex]
\raisebox{1.5ex}{Example 2} & $\beta_{2}$ for $p$ & $2.3462$ & $3.4302$  \\ [0.5ex] \hline \hline 
\end{tabular}
\caption{Order of convergence $\beta_{i},\ (i = 1,2),$ for velocity $\textbf{u}$ and pressure $p$ in Examples $1$ and $2$.}
\label{tab:herrbeta1}
\end{table}
\begin{table}[!ht]
\centering % centering table
\begin{tabular}{ c | c | c | c }
 \hline\hline % inserting double-line
   & Error Estimate & $\textbf{P}_{2}$ & $\textbf{P}_{3}$
 \\ [0.5ex]
\hline    % inserts single-line
% Entering 1st row
& $\beta_{1}$ for $\psi$ & $2.7001$ & $3.7500$ \\ [-0.7ex] 
\raisebox{1.5ex}{Example 1} & $\beta_{1}$ for $c$ & $2.9278$ & $3.8739$  \\ [0.5ex]  \hline 
% inserts single-line
% Entering 1st row
&  $\beta_{1}$ for $\psi$ & $2.9200$ & $3.8200$ \\ [-0.7ex]
\raisebox{1.5ex}{Example 2} & $\beta_{1}$ for $c$ & $2.8972$ & $4.0681$  \\ [0.5ex] \hline \hline 
\end{tabular}
\caption{Order of convergence $\beta_{1}$ for phase-field $\psi$ and concentration $c$ in Examples $1$ and $2$.}
\label{tab:herrbeta2}
\end{table}

\begin{table}[!ht]
\centering % centering table
\begin{tabular}{ c | c | c | c }
 \hline\hline % inserting double-line
   & Error Estimate & $\textbf{P}_{2}-\textbf{P}_{1}$ & $\textbf{P}_{3}-\textbf{P}_{2}$
 \\ [0.5ex]
\hline    % inserts single-line
% Entering 1st row
& $\alpha$ for $\textbf{u}$ & $1.1494$ & $1.1446$ \\ [-0.7ex] 
\raisebox{1.5ex}{Example 1} & $\alpha$ for $p$ & $1.0733$ & $1.0634$  \\ [0.5ex]  \hline 
% inserts single-line
% Entering 1st row
&  $\alpha$ for $\textbf{u}$ & $0.9011$ & $0.9152$ \\ [-0.7ex]
\raisebox{1.5ex}{Example 2} & $\alpha$ for $p$ & $1.0032$ & $0.9944$  \\ [0.5ex] \hline \hline 
\end{tabular}
\caption{Order of convergence $\alpha$ for velocity $\textbf{u}$ and pressure $p$ in Examples $1$ and $2$.}
\label{tab:tauerralpha1}
\end{table}
\begin{table}[!ht]
\centering % centering table
\begin{tabular}{ c | c | c | c }
 \hline\hline % inserting double-line
   & Error Estimate & $\textbf{P}_{2}$ & $\textbf{P}_{3}$
 \\ [0.5ex]
\hline    % inserts single-line
% Entering 1st row
& $\alpha$ for $\psi$ & $1.0558$ & $1.0396$ \\ [-0.7ex] 
\raisebox{1.5ex}{Example 1} & $\alpha$ for $c$ & $1.0602$ & $1.0565$  \\ [0.5ex]  \hline 
% inserts single-line
% Entering 1st row
&  $\alpha$ for $\psi$ & $0.9821$ & $0.9856$ \\ [-0.7ex]
\raisebox{1.5ex}{Example 2} & $\alpha$ for $c$ & $0.9792$ & $0.9815$  \\ [0.5ex] \hline \hline 
\end{tabular}
\caption{Order of convergence $\alpha$ for phase-field $\psi$ and concentration $c$ in Examples $1$ and $2$.}
\label{tab:tauerralpha2}
\end{table}
\subsection{Numerical stability analysis}\label{stabAna}
\begin{figure}[!ht]
\begin{center}
\includegraphics[width=6.5cm,height=5.5cm]{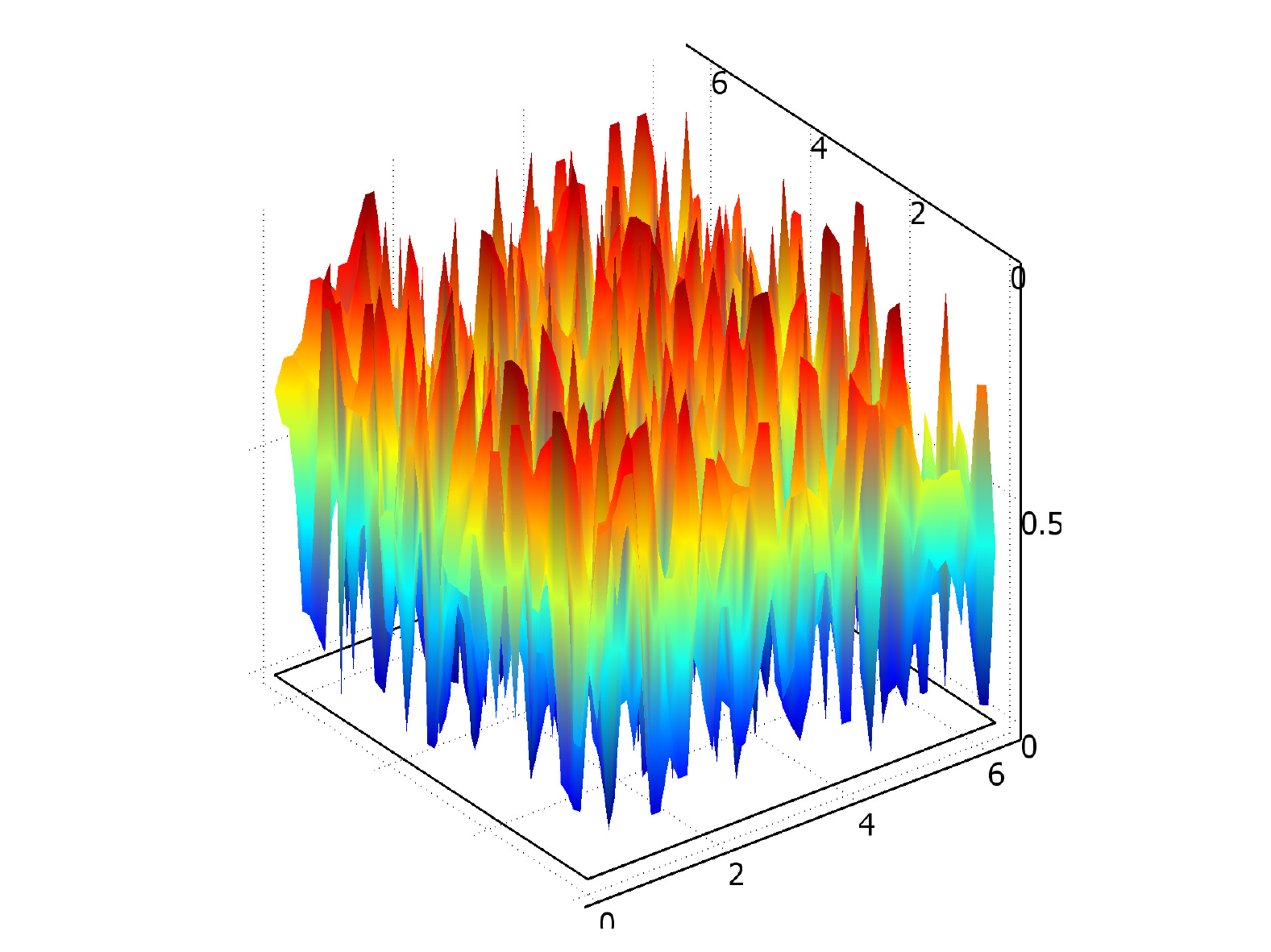} 
\caption{Random function.}
\label{randfn}
\end{center}
\end{figure}

\begin{figure}[!ht]
\centering
\subfigure[Error curves in Example 1]{
   \includegraphics[width=5.5cm,height=5.5cm]{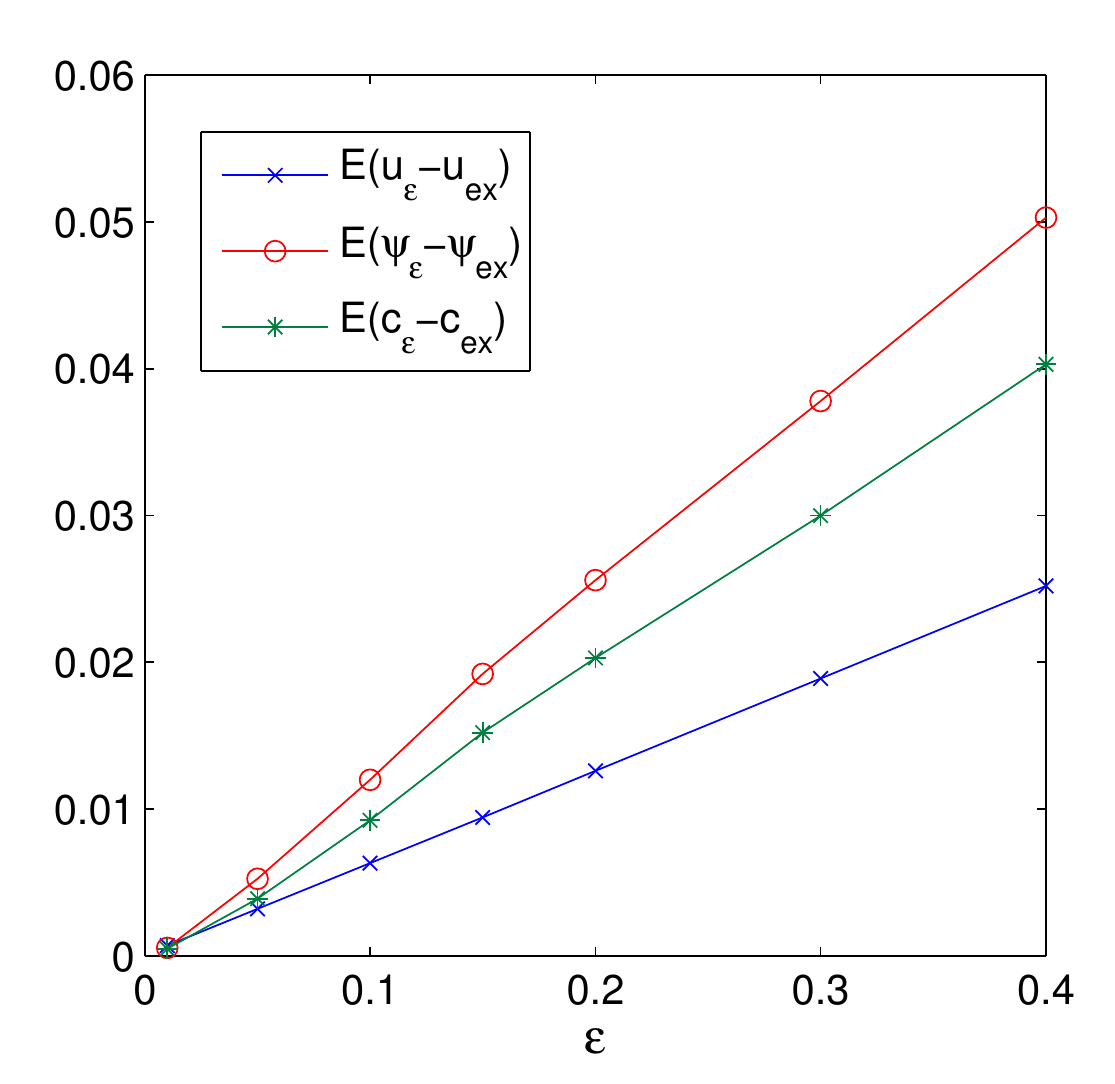}
   \includegraphics[width=5.5cm,height=5.3cm]{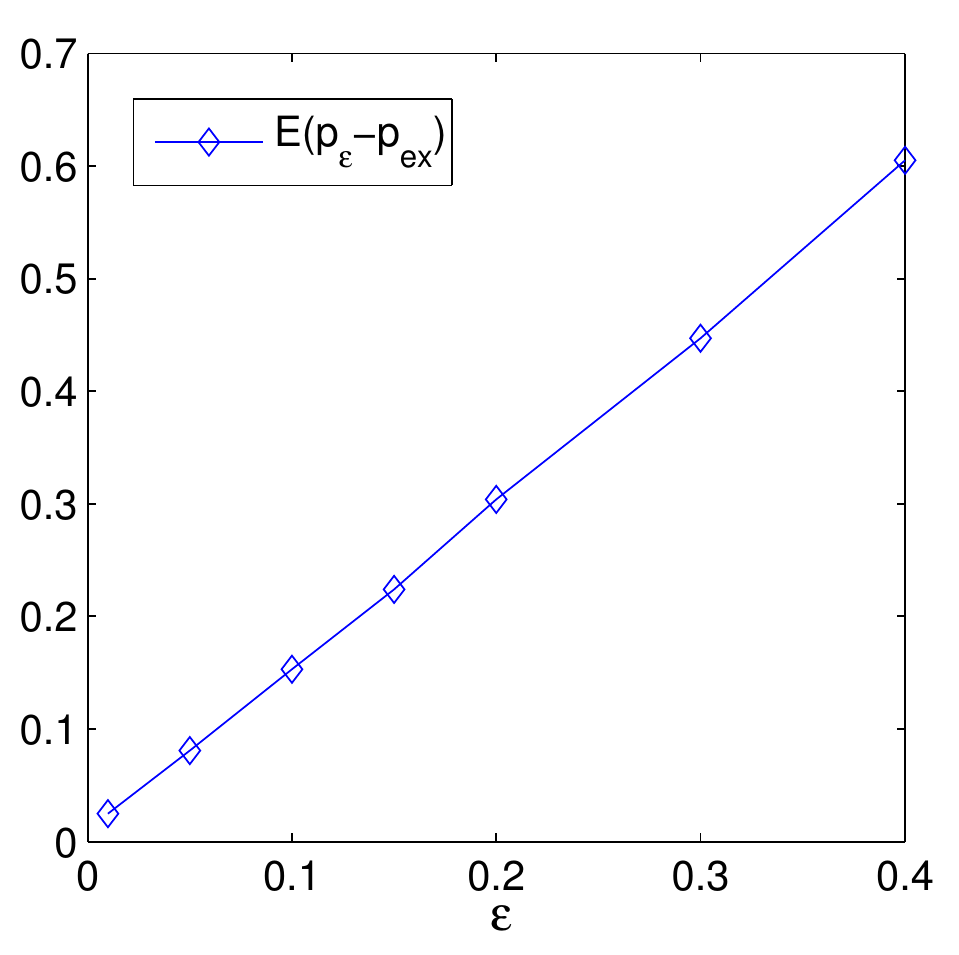} 
   \label{fig:E(Ueps-Uex)a}}
 \subfigure[Error curves in Example 2]{ 		 \includegraphics[width=5.5cm,height=5.5cm]{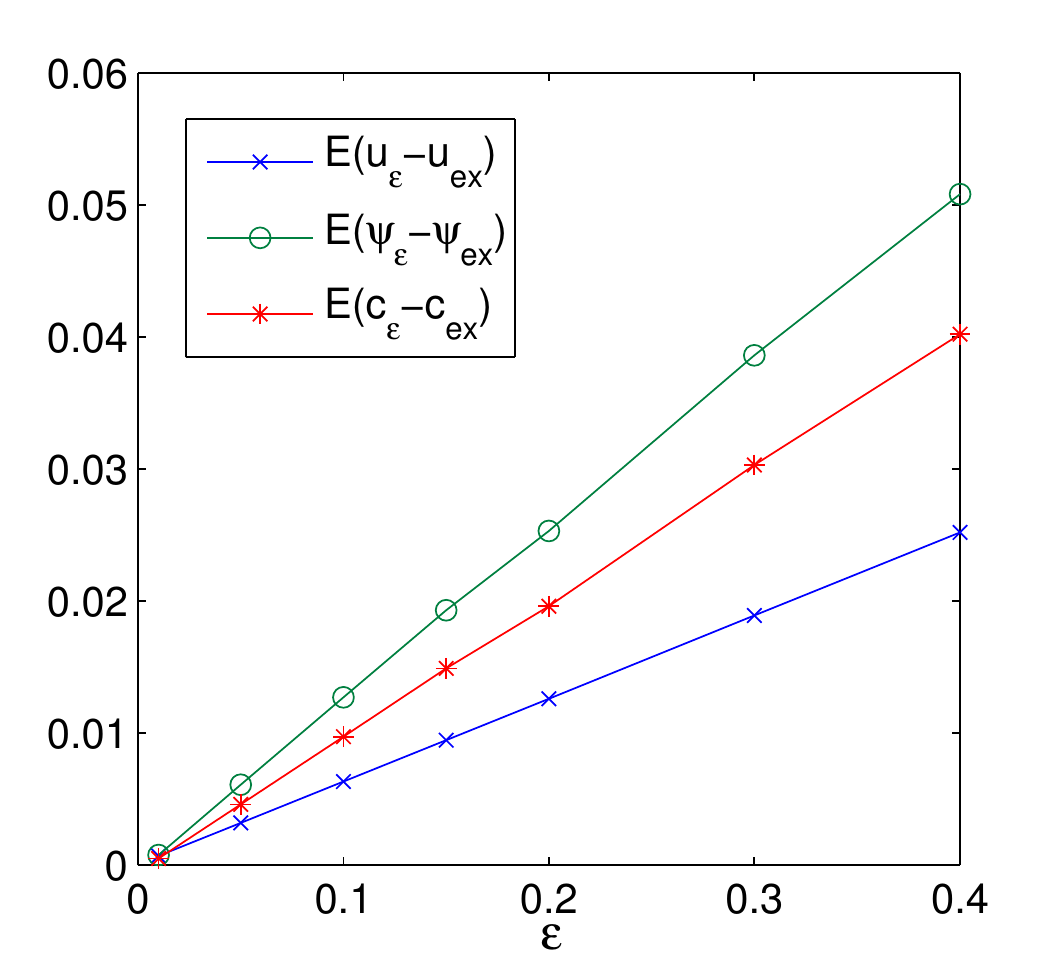}
   \includegraphics[width=5.5cm,height=5.35cm]{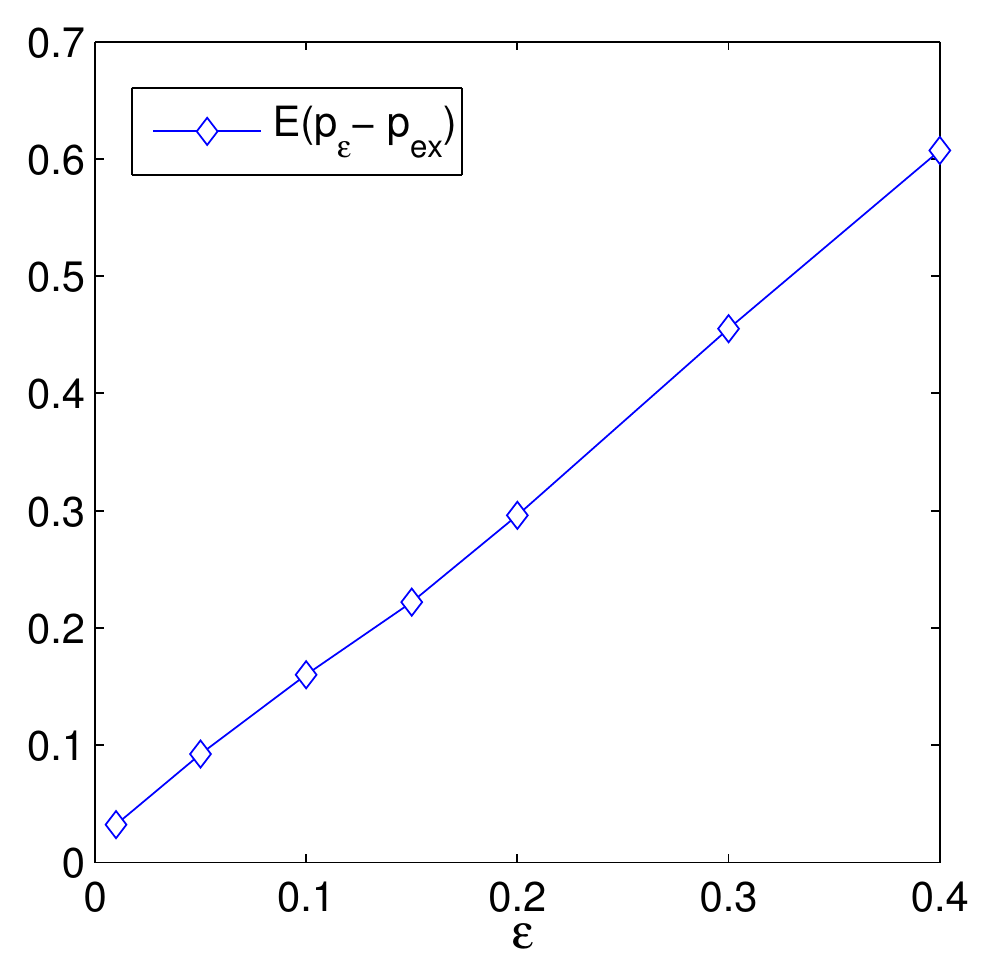} 
   \label{fig:E(Ueps-Uex)b}}
\caption{Errors Curves of norm $EE_{\epsilon,ex}$.}\label{fig:E(Ueps-Uex)}
\end{figure}

\begin{figure}[!ht]
\centering
\subfigure[Error curves in Example 1]{
   \includegraphics[width=5.5cm,height=5.5cm]{Ex1_Uepsil-Uex}
   \includegraphics[width=5.5cm,height=5.4cm]{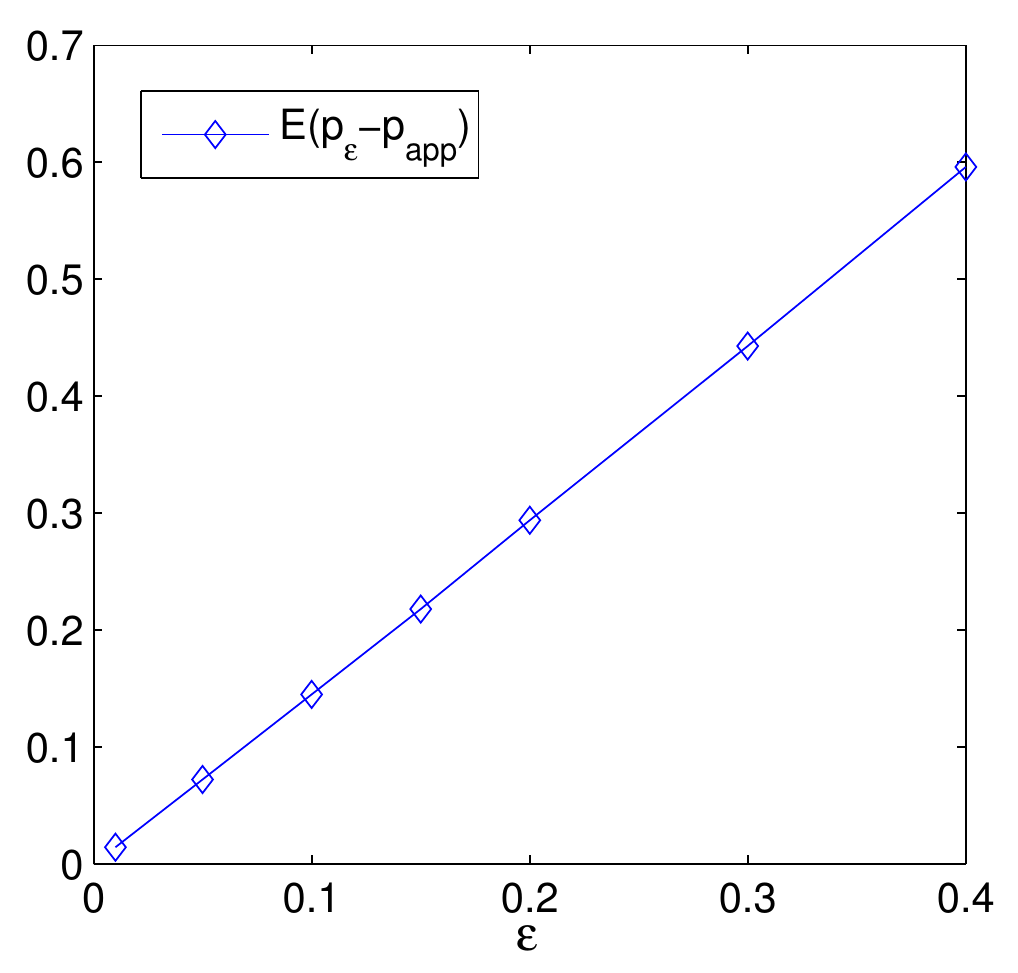} 
   \label{fig:E(Ueps-Uapp)a}}
 \subfigure[Error curves in Example 2]{ \includegraphics[width=5.5cm,height=5.5cm]{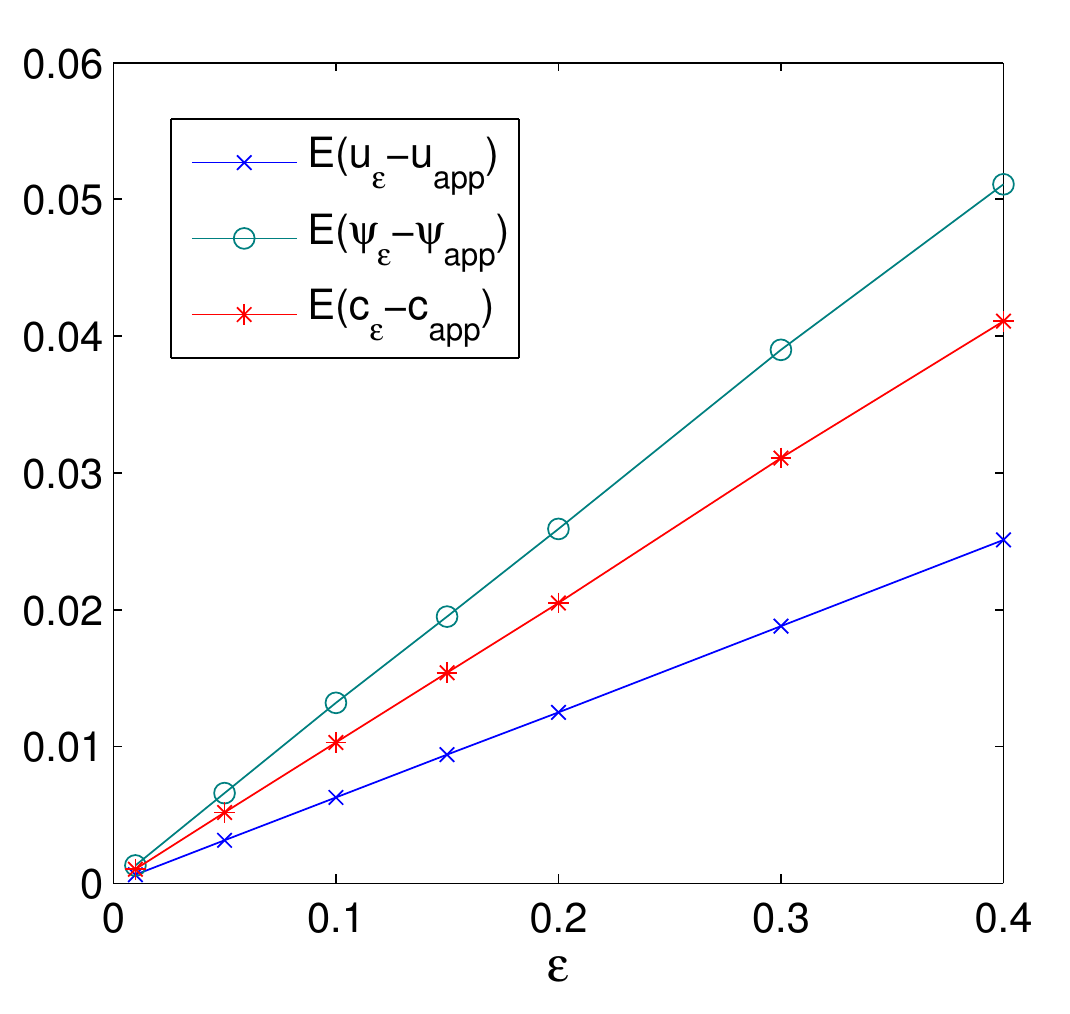}
   \includegraphics[width=5.4cm,height=5.4cm]{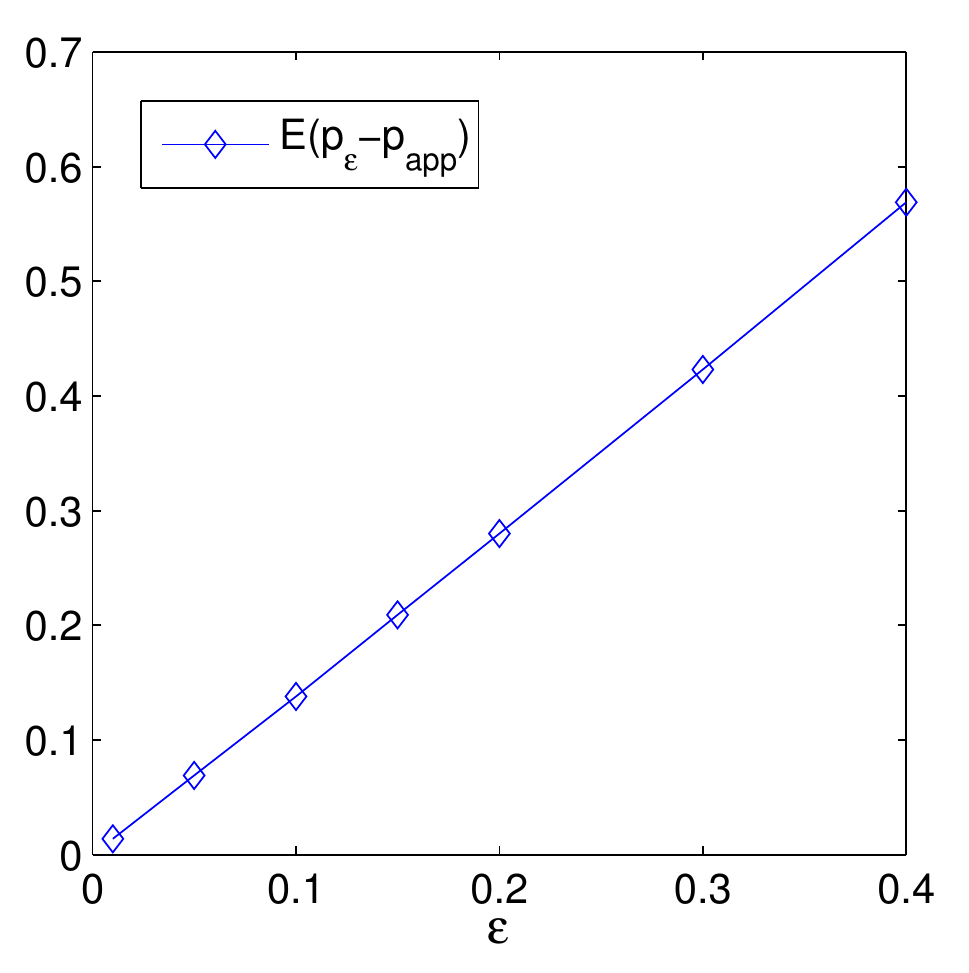} 
   \label{fig:E(Uepsil-Uapp)b}}
\caption{Errors Curves of norm $E_{\epsilon,app}$.}\label{fig:E(Ueps-Uapp)}
\end{figure}

\begin{figure}[!ht]
\begin{center}
\includegraphics[width=7.5cm,height=5.5cm]{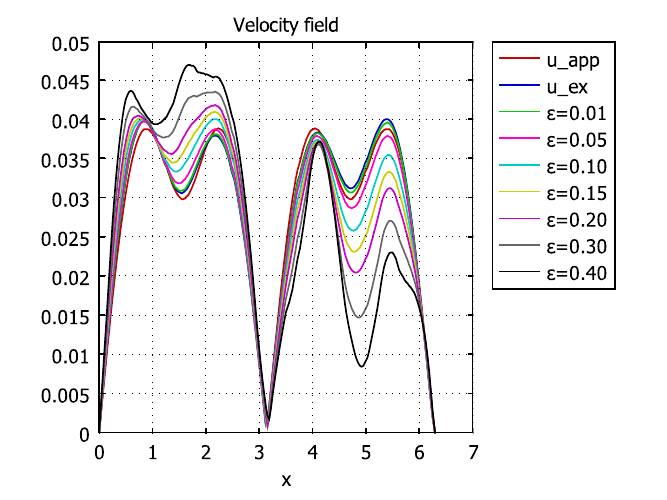}\hspace{-0.35cm}
\includegraphics[width=7cm,height=5.5cm]{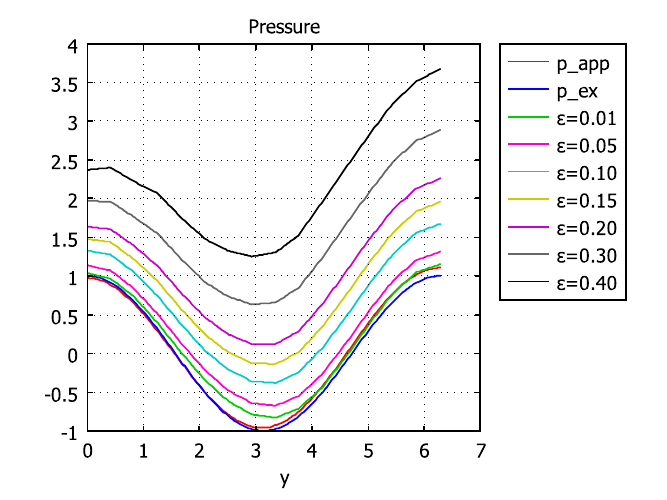}
\includegraphics[width=7cm,height=5.5cm]{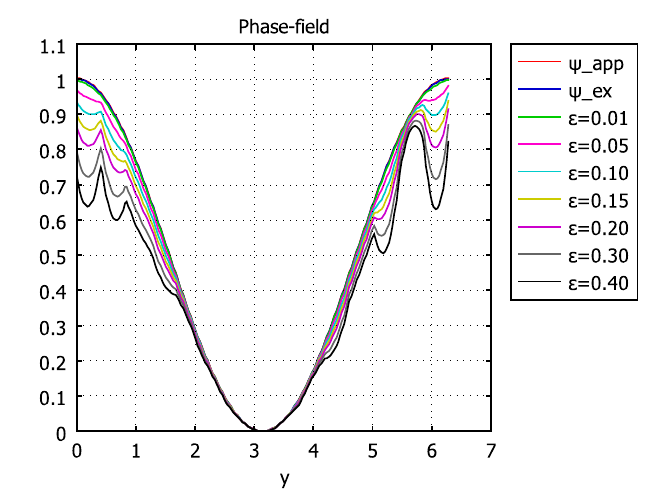}
\includegraphics[width=7cm,height=5.5cm]{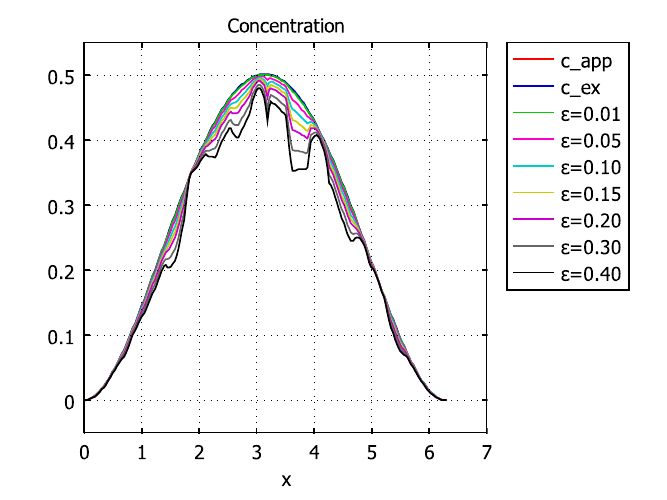}
\caption{Solution curves for the different values of $\epsilon$ in Example $2$.}
\label{Fig:Ex1MultiSolutions}
\end{center}
\end{figure}

\begin{figure}[!ht]
\begin{center}
\includegraphics[width=7cm,height=5.5cm]{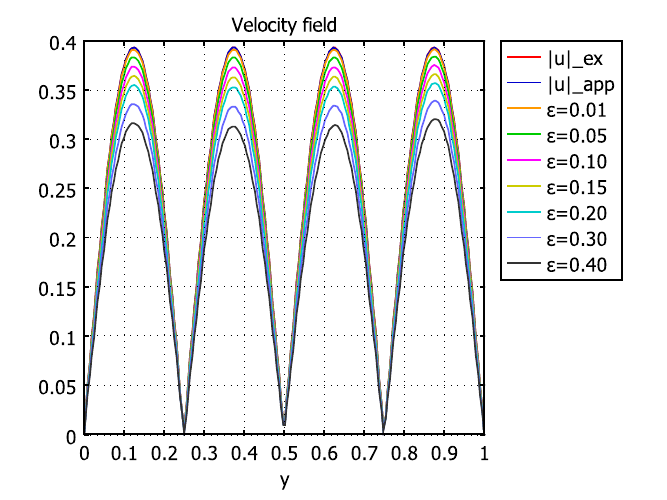}
\includegraphics[width=7cm,height=5.5cm]{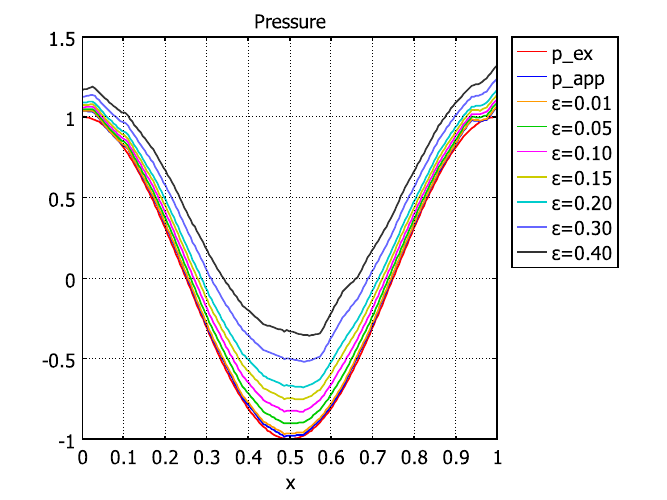}
\includegraphics[width=7cm,height=5.5cm]{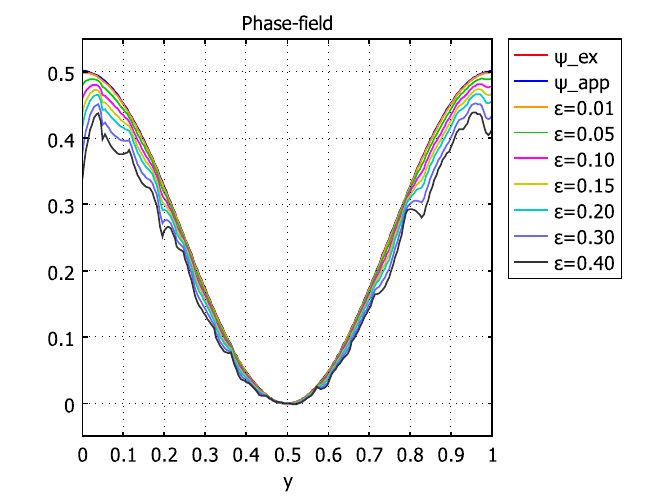}
\includegraphics[width=7cm,height=5.5cm]{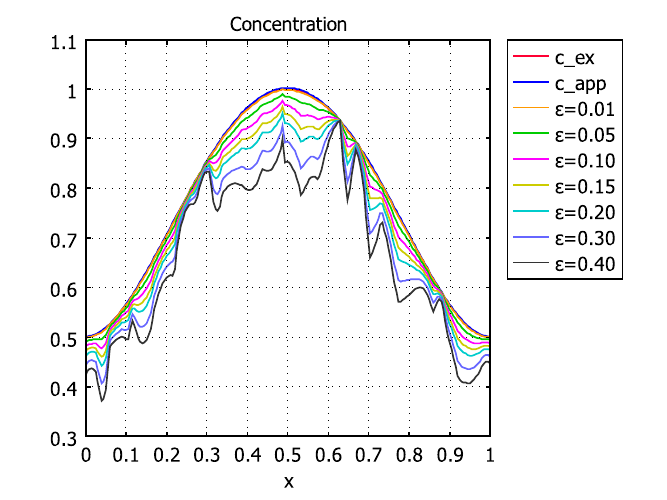}
\caption{Solution curves for the different values of $\epsilon$ in Example $1$.}
\label{Fig:Ex1StabEpsilonRandom}
\end{center}
\end{figure}

An extra term  $(1-\epsilon\ randf)$ has been incorporated in the right hand side terms of each equation in (\ref{FNumScheme}) to introduce $\epsilon$-perturbations, 
where $randf$ assumes values in $[0,1]$ (see e.g. Fig. \ref{randfn}) and $\epsilon$ is the perturbation control parameter. We fix $h=0.2$ and  $\tau=0.1$ and we use quadratic finite elements $\mathbb{P}_{2}$ for $\psi$, $c$ and $\textbf{u}$ and linear finite elements $\mathbb{P}_{1}$ for $p$ in order to study the stability of the method. We perform three different computational stability tests 
(in Fig.\ref{fig:E(Ueps-Uex)},  Fig.\ref{fig:E(Ueps-Uapp)}, Fig.\ref{Fig:Ex1MultiSolutions} and Fig.\ref{Fig:Ex1StabEpsilonRandom}).

In Fig. \ref{fig:E(Ueps-Uex)}, the $L^{2}({\cal Q})$-norm of the discrepancy between exact solution $\Phi_{ex}=(\Phi_{ex}^{(s)})_{s=1,4}=(\textbf{u}_{ex},p_{ex},\psi_{ex},c_{ex})$ and its $\epsilon$-perturbation $\Phi_{\epsilon}=(\Phi_{\epsilon}^{(s)})_{s=1,4}=(\textbf{u}_{\epsilon},p_{\epsilon},\psi_{\epsilon},c_{\epsilon})$ i.e. $E_{\epsilon,ex}\left(\Phi_{\epsilon}^{(s)}-\Phi_{ex}^{(s)}\right)=\left\|\Phi_{\epsilon}^{(s)}-\Phi_{ex}^{(s)}\right\|_{L_{2}(Q)}$, for $s=1,4$ are plotted versus $\epsilon$ (which are shown for $\epsilon=0.01,0.05,0.1,0.15,0.2,0.3,0.4$). A linear dependence of errors on $\epsilon$ is observed, indeed, $E_{\epsilon,ex}\left(\Phi_{\epsilon}^{(s)}-\Phi_{ex}^{(s)}\right)\approx\ m_{s}\ \epsilon$, for $s=1,4$, where $m_{s}$ represents the slope of the error curve; refer to Table \ref{tab:Staberr}.
In Fig. \ref{fig:E(Ueps-Uapp)}, the error $E_{\epsilon,app}\left(\Phi_{\epsilon}^{(s)}-\Phi_{app}^{(s)}\right)=\left\|\Phi_{\epsilon}^{(s)}-\Phi_{app}^{(s)}\right\|_{L_{2}(Q)}$, for $s=1,4$ between the approximate solution $\Phi_{app}=(\Phi_{app}^{(s)})_{s=1,4}=(\textbf{u}_{app},p_{app},\psi_{app},c_{app})$ without random (i.e., $\epsilon =0)$ and $\Phi_{\epsilon}$, are plotted against $\epsilon$. The same observation holds as
in Fig. \ref{fig:E(Ueps-Uex)}; refer also to Table \ref{tab:Staberr}. 

Finally, in   Fig. \ref{Fig:Ex1MultiSolutions} and  Fig. \ref{Fig:Ex1StabEpsilonRandom}, the solutions are shown on a part of domain and at time $t=1$ in order demonstrate stability with respect to perturbations. In   Fig. \ref{Fig:Ex1MultiSolutions}, we  fix $y = \pi/2$ and $x$ varies for velocity and concentration, and $t = 1$, $x = \pi$ and $y$ varies for pressure and phase field. In   Fig. \ref{Fig:Ex1StabEpsilonRandom} we fix $x = 1/2$ and $y$ varies for velocity and phase field, and $x = 1/2$ and $y$ varies. The graphs shows that the solution is stable.
\begin{table}[!ht]
\centering % centering table
\begin{tabular}{ c | c | c | c }
 \hline\hline % inserting double-line
   & Slope & $E_{\epsilon,ex}$ & $ E_{\epsilon,app}$
 \\ 
\hline \hline    % inserts single-line
% Entering 1st row
& $m_{\textbf{u}}$ & $0.1701$ & $0.1754$ \\ 
& $m_{\psi}$ & $0.8638$ & $0.8818$ \\ 
\raisebox{0.5ex}{Example 1} & $m_{c}$ & $0.4341$ & $0.4371$  \\ & $m_{p}$ & $1.4738$ & $1.4913$ \\  \hline 
% inserts single-line
% Entering 1st row
&  $m_{\textbf{u}}$ & $0.0628$ & $0.0635$ \\ 
&  $m_{\psi}$ & $0.1283$ & $0.1347$ \\ 
\raisebox{0.5ex}{Example 2} & $m_{c}$ & $0.1018$ & $0.1065$  \\  
&  $m_{p}$ & $1.4877$ & $1.4236$ \\  \hline \hline 
\end{tabular}
\caption{Slopes of Norm $L_{2}$ in Examples $1$ and $2$.}
\label{tab:Staberr}
\end{table}
%____________________________________________________________________
%
%              Example Anisotropic Multiple Solution Curves
%_____________________________________________________________________
%
\section{Concluding remarks}\label{discussions}
This paper presents a numerical investigation and resolution of the isothermal anisotropic solidification model (\ref{MP}). The purpose of this study is to validate the derived numerical scheme by performing its error and stability analysis. The model has been discretized  with respect to spatial and time variables . The discretization result into a system of nonlinear ordinary differential equations. The resulting non-linear systems are solved by using a solver DASSL. 

Second, the convergence and stability of the numerical scheme has been validated (both with respect to space and time variables) by considering two examples with known exact solutions. It is numerically  demonstrated that the error estimates with respect to spatial coordinates are of order $i+1$ for ${\bf u}$,  $\psi$  and  $c$ and of order $i$ for $p$, and the error estimates for time are of order 1 for $({\bf u},  p,\psi, c)$. Both of these numerical estimates are in excellent accordance with the estimates (\ref{CF1}). The stability of the scheme has also been verified by introducing a random function, which varies between 0 and 1, in the model. It is found that the numerical scheme is completely stable and it has linear dependence with the increase in percentage error.

\section{Declarations}
\textbf{Competing interests:} The author declare that he has no competing interests.\newline 
\textbf{Acknowledgments:} The author is thankful to the journal for consideration of this manuscript.

\end{document}